\documentclass[11pt]{article}
\usepackage{graphicx,amsmath,amssymb,mathrsfs}
\usepackage[a4paper,  margin=1in]{geometry}
\usepackage[utf8]{inputenc}
\usepackage[english]{babel}
\usepackage[utf8]{inputenc}
\usepackage{amsmath}
\usepackage{longtable}
\usepackage{amsfonts}
\usepackage{amssymb} 
\usepackage{enumerate}
\usepackage{enumitem} 
\usepackage{afterpage}
\usepackage{amsmath}
\usepackage{float}
\usepackage{amsthm}
\usepackage{enumerate}
\usepackage{enumitem}
\usepackage{pifont}
\usepackage{amssymb}
\usepackage{emptypage}
\usepackage{bigints}
\usepackage [autostyle, english = american]{csquotes}
\MakeOuterQuote{"}

\usepackage{graphicx}
\usepackage{rotating}
\usepackage{wrapfig}
\usepackage{caption}
\usepackage{stackengine}
\usepackage{subcaption}
\usepackage{xspace}
\usepackage{graphicx,amsmath,amssymb,mathrsfs,tcolorbox}
\usepackage[a4paper,  margin=1in]{geometry}
\usepackage{authblk}
\usepackage{natbib}
\usepackage[none]{hyphenat}
\usepackage{setspace}
\setcitestyle{authoryear,open={(},close={)}}

\newtheorem{theorem}{Theorem}[section]

\usepackage{hyperref}
\hypersetup{
    colorlinks=true,
    linkcolor=blue,
    filecolor=magenta,      
    urlcolor=cyan,
    citecolor=blue
    }

\begin{document}

\title{\textbf{Estimation of Proportion of Null Hypotheses Under Dependence}}

\author{Nabaneet Das}
\affil{Indian Statistical Institute,Kolkata}

\maketitle

\begin{abstract}
Estimation of the proportion of null hypothesis ($\pi_0$) in a multiple testing problem can greatly enhance the performance of the existing algorithms. Although various estimators for $\pi_0$ have been proposed, most are designed for independent samples, and their effectiveness in dependent scenarios is not well explored. This article investigates the asymptotic behavior of the BH estimator and evaluates its performance across different types of dependence. Additionally, we assess Storey's estimator and another estimator proposed by \cite{patra2016estimation} to understand their effectiveness in these settings.

\end{abstract}

\maketitle

\section{Introduction}
It has been shown that, the FDR controlling procedure of Benjamini-Hochberg algorithm controls FDR at a level less than or equal to $ \frac{n_0}{n} \alpha = \pi_0 \alpha $ under independence (\cite{benjamini1995controlling}) (Here $\pi_0 $ represents the proportion of true null hypotheses). Traditional methods often reach their maximum Type I error rate when faced with the global null hypothesis, leaving ample room for increased statistical power, particularly when $\pi_0 < 1$. In an ideal scenario where $\pi_0$ is known, adjusting the Benjamini-Hochberg procedure to precisely control the FDR at a specified $\alpha$ level would be straightforward. However, in practice, $\pi_0$ remains unknown and necessitates estimation from the available data. \cite{benjamini2000adaptive} proposed an adaptive FDR controlling procedure where $\pi_0$ was estimated based on a graphical approach similar to \cite{schweder1982plots}. Building upon this foundation, \cite{storey2002direct} later proposed an adaptive technique grounded in a conservative estimation of $\pi_0$, showcasing substantial enhancements in statistical power. A bias-corrected version of Storey's estimator was later proposed by \cite{biswas2022bias}. \cite{genovese2004stochastic} brought attention to the identifiability challenges inherent in the estimation of $\pi_0$, delineating key criteria for achieving identifiability and purity in estimations. In the context of their study, a comprehensive review of existing literature was provided, synthesizing the advancements in FDR control methodologies. Furthermore, a novel estimator for $\pi_0$ was proposed, accompanied by a robust confidence interval construction. Notably, \cite{genovese2004stochastic} demonstrated the asymptotic validity of both the proposed point estimator and associated confidence intervals, reinforcing their utility in practical applications. \cite{efron2012large} has addressed the identifiability issue by assuming that the support of the null density contains a portion which is disjoint from the non-null density. Many works have been done under the assumption that the distributions belong to certain parametric models (see, for example \cite{cohen1967estimation}, \cite{lindsay1993multivariate}). \cite{jin2008proportion}, \cite{cai2010optimal} used empirical characteristic functions to estimate the non-null distribution under a semiparametric normal mixture model. In the context of multiple testing \cite{meinshausen2005lower}, \cite{meinshausen2006estimating}, \cite{celisse2010cross} have also suggested different estimators and confidence bounds. \cite{patra2016estimation} have investigated the identifiability problem in complete generality and proposed a consistent estimator of $\pi_0$. \\ 
 Although these estimators were initially conceived under the assumption of independence among test statistics (such as p-values) corresponding to each hypothesis, scant literature exists on exploring the extent to which these estimators can withstand dependence. In this chapter, we delve into an investigation of the consistency of these estimators across a spectrum of weakly dependent scenarios. 
 \section{The model and identifiability}
 We consider a sample $X_1, X_2, \dots , X_n $ which comes from a common distribution function $F$. In the multiple testing problem with $n$ null hypotheses $H_{01}, \dots , H_{0n} $, $X_i$ is the test statistic corresponding to the $i$-th testing problem. Let $F_0$ and $F_1$ be the distribution under the null hypothesis and alternate hypothesis respectively ($F_0$ and $F_1$ does not depend on $i$). Then we consider the mixture model 
 \begin{equation}\label{mixture}
     F(x) = (1-\pi_1) F_0 (x) + \pi_1 F_1 (x) 
 \end{equation}
 Note that, the mixture model is written as a function of proportion of non-null hypotheses ($ \pi_1 = 1 - \pi_0$).
 Usually the null distribution $F_0(.) $ is known and $\pi_1 , F_1 (.) $ are unknown. For example, if we consider p values for the multiple testing problem, then $F_0 $ is the $Unif [0,1] $ distribution if the p values are continuous. If we consider z values as per \cite{efron2008microarrays}, then the null distribution is standard normal distribution. \\
 When both $\pi_1$ and $F_1$ are unknown, the problem is non-identifiable. If the model \ref{mixture} holds then for any $ 0 \leq \gamma \leq \pi_0 = 1 - \pi_1 $, one can re-write as 
 $$ F(x) = (\pi_1 + \gamma ) \Bigg( \left( \frac{ \pi_1 }{ \pi_1 + \gamma } F_1 (x) + \frac{\gamma}{ \pi_1 + \gamma  } F_0 (x) \right) + ( 1 - \pi_1 - \gamma) F_0 (x) \Bigg) $$
 So, the expression of $F(x)$ as a convex combination of the null distribution ($F_0(x) $) and some other distribution is not uniquely defined. \cite{genovese2004stochastic} discussed the identifiability issue when $F_0$ is $Unif[0,1] $ distribution and $F$ has a density. Building on this foundation, \cite{patra2016estimation} redefine the mixing proportion as follows 
 \begin{equation}\label{proportion}
  \text{  \underbar{$\pi_1$}}  = \inf \: \{ \: \gamma \in (0, 1 ] \: : \: \Big\{ \frac{ F - (1- \gamma) F_0 }{ \gamma } \Big\} \text{ is a CDF} \} 
 \end{equation}
They have further examined this $\text{  \underbar{$\pi_1$}}$ and found the connection between $\text{  \underbar{$\pi_1$}}$ and $\pi_1$ satisfying model \ref{mixture}. 
 \begin{tcolorbox}

 \begin{theorem}\label{Expression}{\cite{patra2016estimation} }

Let F be as in model \ref{mixture} and $\text{  \underbar{$\pi_1$}}$ be defined as \ref{proportion}. Then, 
\begin{equation}\label{sub-cdf}
\text{  \underbar{$\pi_1$}}= \pi_1 - \sup \{ \: 0 \leq \epsilon \leq 1 \: : \: \pi_1 F_1 - \epsilon F_0 \text{ is a sub-CDF } \}  
\end{equation}
In particular, if $F_1$ and $F_0$ are absolutely continuous, i.e. they have densities $f_1$ and $f_0$ respectively, then 
\begin{equation}\label{essinf}
\text{  \underbar{$\pi_1$}} = \pi_1 \Big(  1 - \text{ essinf }\frac{f_1}{f_0}\Big)   
\end{equation}
where, for any function $ g , \text{essinf } g = \sup \{ a \in \mathbb{R} \: : \: m \left\{  \{ x \: : \: g(x) < a  \right\} = 0 \} $, $m$ being the Lebesgue measure. 
\end{theorem}
\end{tcolorbox} 
As a consequence, $\text{  \underbar{$\pi_1$}} < \pi_1 $ if and only if there exists $ c > 0 $ such that $ f_1 \geq c f_0 $, almost everywhere $m$. Note that, in particular, if the support of $F_1$ is strictly contained in that of $F_0$, then the problem is identifiable and we can estimate $\pi_1$. In the subsequent sections we will discuss different estimators of $\pi_1$ and their consistency.  
\section{ Inconsistency of the BH estimator}
The Benjamini Hochberg step up algorithm is based on ordered p values $P_{ (1) } \leq \dots \leq P_{ ( n ) } $ which compares each $P_{ (i) } $ to $ \frac{ i \alpha}{n} $. When the observations are independent, this algorithm controls the FDR at a level $ \pi_0 \alpha \leq \alpha$. Later an adaptive version of BH procedure was proposed by \cite{benjamini2000adaptive} and they proposed a graphical approach of estimating $n_0$ and hence $ \pi_0 = \frac{n_0}{n}$. The estimation works as follows 
\begin{itemize}
    \item Calculate the slopes $ S_i = \frac{ (1 - P_{ (i) } ) }{ (n+1 -i) } $ 
    \item Let, $ j = \inf \{ i \geq 2 \: : \: S_i < S_{i-1} \}$. 
    Then set, $\hat{n}_0 = \min ( [\frac{1}{S_j} +1 ],n) $ and $ \hat{\pi}_0  = \frac{\hat{n}_0}{n}$
\end{itemize}
One might expect that the estimator $\hat{\pi}_{0,BH} $ is a consistent estimator of $\pi_0$ in order to ensure a good performance of the adaptive BH method for large number of hypotheses. However, the following theorem describes that $\hat{\pi}_{0,BH} $ cannot be a consistent estimator of $\pi_0$ as it puts non-zero mass on some point for any $0 \leq \pi_0 \leq 1 $. 

\begin{theorem}\label{inconsistent_BH}
 Let $j$, $\hat{n}_0 $ and $\hat{\pi}_{0,BH}  $ be defined  as per the above mentioned description. Also assume that, the p values are absolutely continuous and the alternate distribution of the p values are stochastically smaller than the null distribution ($Unif[0,1]$). Then, $ \hat{\pi}_{0,BH}  \nrightarrow \pi_0 $ for all values of $ \pi_0 \in [ 0, 1 ]$ (i.e. $ \hat{\pi}_{0,BH}  $ is not a consistent estimator of $ \pi_0$) 
\end{theorem} 
\subsection{Proof of theorem \ref{inconsistent_BH}} 
By definition, $ j \geq 2 $ with probability 1. \\[0.1 in]
$
P( j > 2 ) \\
= P( S_2 \geq S_1) \\[0.05 in]
= P( \frac{ 1 - p_{ (2) }}{ n -1 } \geq \frac{ 1 - p_{ (1) }}{n} ) 
$ \\[0.05 in]
Let, $q_i = 1 - p_{ i} $ for $ i = 1 , \dots , n $. This implies, $ q_{ (n + 1 - i) } = 1 - p_{ (i) } $ for all i. \\[0.05 in]
So, $ P( S_2 \geq S_1 ) = P( \frac{ q_{ (n) }}{ q_{ (n-1) }} \leq \frac{n}{n-1 })$ \\[0.05 in]
Let $g$ and $G$ respectively denote the density and the cdf of the $q_i$'s. Then, from the joint density of order statistics, one can easily deduce that, \\[0.1 in]
$ P( q_{ (n-1)  } \geq ( \frac{n-1}{n} ) q_{ (n) } ) \\[0.05 in]
= n ( n-1) \bigint\limits_{0}^{1}  \Big(  \int\limits_{ (\frac{n-1}{n})y}^{ y} (G(x) )^{n-2} g(x) dx  \Big) g(y) dy \\[0.05 in ]
= \bigint\limits_{0}^{1} n \Big( (G(y))^{n-1} - (G \left( ( \frac{n-1}{n} ) y \right) )^{ n-1} \Big) g(y) dy \\ [0.05 in] 
= \bigint\limits_{0}^{1} \Big( 1 - \left( \frac{ G( ( \frac{ n-1}{n} ) y )}{ G(y) } \right)^{ n-1} \Big) \left( n ( G(y))^{n-1} g(y) \right) dy \\[0.05 in]
= E \Big[  1 - \left(  \frac{ G ( (\frac{ n-1}{n }) q_{ (n ) } )}{ G(q_{ (n) })} \right)^{ n-1 } \Big]
$ \\[0.2 in]
The last equality holds because of the fact that $q_{ (n) } $ has density $ h(y) = n \left( G(y) \right)^{n-1} g(y) I_{ 0 < y < 1 } $. Now observe that, \\[0.1 in] 
$
n \Big[ 1 - \frac{ G \left( (\frac{n-1}{n}) q_{ (n)  } \right) }{ G \left( q_{ (n) }  \right)  }\Big] \\[0.05 in]
= n \Bigg[ \frac{ G \left( q_{ (n) } \right) - G \left( (\frac{n-1}{n})q_{ (n) } \right) }{ G \left( q_{ (n) } \right) } \Bigg] \\[0.1 in]
= \frac{ q_{ (n) } g ( \xi_n) }{ G \left( q_{ (n) } \right)} 
$  for some $(\frac{ n-1}{n} )q_{ (n) } < \xi_n < q_{ (n) }$\\[0.05 in]
Since the alternate density is stochastically smaller than the $Unif [0,1] $ distribution, $ G(1) =1 $ and $ q_{ (n) } \stackrel{a.s.}{\to} 1 $ as $ n \to \infty$. \\[0.05 in]
So, by the continuity of $G(.) $, we can write $ G(q_{ (n) }) \stackrel{a.s.}{\to}  G(1) = 1 $ as $ n \to \infty$. \\
And $ g ( \xi_n) \to g(1-) $ as $ n \to \infty $. \\
Thus, 
$$ n \Big[ 1 - \frac{ G \left( (\frac{n-1}{n}) q_{ (n)  } \right) }{ G \left( q_{ (n) }  \right)  }\Big] \to g(1-) \text{  as } n \to \infty  $$
Here, $ g(1-) = \lim_{ x \uparrow 1  } g(x) $. \\ [0.1 in]
And hence, $$ \lim_{ n \to \infty} P( j > 2 ) = \lim_{ n \to \infty} E \Big[  1 - \left(  \frac{ G ( (\frac{ n-1}{n }) q_{ (n ) } )}{ G(q_{ (n) })} \right)^{ n-1 } \Big] = 1 - \exp ( -  g( 1-) )  $$ 
\textbf{ \underline{Case 1}} : $g(1-) < \infty $ and $1-\exp(-g(1-)) < 1 $ \\

Since $ P( j \geq 2 ) = 1 $ and $ \lim_{ n \to \infty} P( j > 2 ) < 1 $, this implies, 
$$ \lim_{ n \to \infty} P(j =2) = \exp ( - g(1-) ) > 0 $$ 
And thus, $$ \lim_{ n \to \infty } P( \hat{n}_0 = \frac{1}{ S_2} +1 )   \geq \exp ( - g(1-)) > 0 $$ 
Since $ S_2 = \frac{1 - p_{(2) }}{n-1} $, this implies, 
$$ \lim_{ n \to \infty } P( \hat{n}_0 = n ) > 0 $$
and hence 
$$ \lim_{ n \to \infty } P( \hat{\pi}_{0,BH}  = 1 ) > 0 $$
So, even if the true value of $ \pi_0 $ is strictly less than 1, the estimator asymptotically puts a positive mass at 1 and hence, this estimator is not a consistent estimator for $\pi_0$.  \\
\textbf{ \underline{Case 2}} : $g(1-) = \infty $ and $1-\exp(-g(1-)) = 1 $ \\
For $ i = 0, \dots , n-2 $, define 
$$ A_i =  \Big\{ \frac{ q_{ (n- i)}}{ n-i} \leq \frac{ q_{(n-i-1)}}{n-i-1} \Big\}
$$
Then, $$ P( j > k ) = P( S_k \geq S_{k-1} \geq \dots \geq S_1) = P \left( \bigcap\limits_{i=0}^{k-2} A_i \right) $$
We have already shown that $ \lim_{ n \to \infty} P(A_0) = 1 - \exp (- g(1-) ) = 1 $.\\ 
\medskip
From the the joint density of $ (q_{ n-i} , q_{ n-i-1} )$, it can be deduced that, \\[0.05 in] 
$P \left( \frac{ q_{(n-i) }}{n-i} \leq \frac{q_{ (n-i-1}}{n-i-1} \right)$ \\[0.05 in]
$ = \frac{ n!}{ (n-i-2)! i! } \bigint\limits_{0}^{1} \bigint\limits_{ \left( \frac{n-i-1}{n-i} \right) y }^{y} \left( G(x) \right)^{i-2} \left( 1 - G(y) \right)^{ i} g(x) g(y) dx dy  $ \\[0.05 in]
$ = \frac{ n!}{ (n-i-1)! i! } \bigint\limits_{0}^{1} \left( 1 - G(y) \right)^i \Big( \left( G(y) \right)^{ n-i-1} - \left( G \left( ( \frac{n-i-1}{n-i} y\right) \right)^{ n-i-1} \Big) g(y) dy $  \\[0.05 in] 
$  = \frac{ n!}{ (n-i-1)! }{i! } \bigint\limits_{0}^{1} \left( G(y) \right)^{ n-i-1} ( 1 - G(y) )^i \Bigg( 1 - \Big( \frac{G \left( \left(  \frac{n-i-1}{n-i} y\right)\right)}{ G(y) } \Big)^{n-i-1} \Bigg) g(y) dy  $ \\[0.05 in]
$ = E \Bigg[ 1 - \Big( \frac{G \left( \frac{n-i-1}{n-i} q_{(n-i)}\right)}{G (q_{ (n-i) })} \Big)^{ n-i-1} \Bigg]$ \\[0.15 in]
Note that, if $i = o(n) $ then $q_{ (n-i) } \stackrel{a.s.}{\to}  1 $ and $ G (q_{ (n-i) }) \stackrel{a.s.}{\to}  1 $ as $ n \to \infty$ . Since $ n \Big( 1 - \frac{G \left( \frac{n-i-1}{n-i} q_{(n-i)}\right)}{G (q_{ (n-i) })} \Big) \stackrel{a.s.}{\to}  g(1-) = \infty$, this implies, $\lim_{ n \to \infty} P(A_i) = 1 $. \\[0.1 in] 
If $ \lim_{ n \to \infty} \frac{i}{n} = \delta \in (0,1) $, then if $z_{ \delta} $ is the $(1- \delta)$-th quantile of the distribution of q's, then $  q_{ (n-i) } \stackrel{a.s.}{\to}  z_{ \delta} $ as $ n \to \infty$. \\[0.05 in] 
Again, by continuity theorem, $ G (q_{ (n-i) } ) \stackrel{a.s.}{\to}  (1 - \delta) $  as $ n \to \infty $ and, 
$$ n \Big( 1 - \frac{G \left( \frac{n-i-1}{n-i} q_{(n-i)}\right)}{G (q_{ (n-i) })} \Big) = \frac{ z_{ \delta} g ( z_{ \delta }-) }{(1- \delta)^2 }$$  
And thus, $$ \lim_{ n \to \infty } E \Bigg[ 1 - \Big( \frac{G \left( \frac{n-i-1}{n-i} q_{(n-i)}\right)}{G (q_{ (n-i) })} \Big)^{ n-i-1} \Bigg] = 1 - \exp \Big( - \frac{ z_{ \delta} g ( z_{ \delta }-) }{(1- \delta)^2 } \Big)$$
Since this is true for all $ 0 < \delta <1 $, we can surely find some $ 0 < \delta_0 < 1 $ for which $ 1 - \exp \Big( - \frac{ z_{ \delta} g ( z_{ \delta }-) }{(1- \delta)^2 } \Big) < 1 $ and for that $ \delta_0$, we'll have 
$$ \lim_{ n \to \infty } P( \hat{\pi}_{0,BH}  = \delta_0) > 0  $$ 
Existence of such $ \delta_0$ will not allow $ \hat{\pi}_{0,BH}  $ to have a degenerate limit and hence it is not a consistent estimator of $ \pi_0$. \\[0.06 in]
Theorem \ref{inconsistent_BH} demonstrates that the BH estimator lacks consistency for $\pi_0$. Despite this limitation, it often serves as a conservative estimator of $\pi_0$ (i.e., $\hat{\pi}_{0,BH}  \geq \pi_0 $). Consequently, the adaptive BH procedure typically maintains validity for FDR control. The following theorem illustrates this fact for the global null. 

\begin{theorem}\label{pi_0=1}
Under the global null (i.e. $\pi_0 = 1 $), the BH estimator of $\pi_0 $ satisfies $$ \hat{\pi}_{0,BH}  \stackrel{a.s.}{\to}  1 \text{  as }n \to \infty  $$  
\end{theorem} 
\textbf{Proof }: If $\pi_0 = 1 $, then all P values come from $Unif [0,1] $ distribution. \\[0.05 in]
In this case, $q_1, \dots , q_n \stackrel{i.i.d.}{\sim}  Unif [0,1] $ where $q_i = 1 - P_i$ for $i=1, \dots , n $. \\[0.05 in]
Observe that, for $ k = 2, \dots , n $ \\[0.05 in] 
$ P( j > k ) $ \\[0.05 in]
$ = P( S_k \geq S_{k-1} \geq \dots \geq S_2 \geq S_1) $ \\[0.05 in]
$ = P \left( \frac{ q_{ (n+1-k) }}{ n+1-k } \geq  \frac{ q_{ (n+2-k) }}{ n+2-k } \geq  \dots \frac{ q_{ (n) }}{ n } \right)$ \\[0.05 in] 
$ = 
n! \bigint\limits_{0}^{1} \bigint\limits_{ \left( \frac{n-1}{n} \right) q_{ (n) }}^{q_{ (n) }} \dots \bigint\limits_{ \left( \frac{n+1-k}{n+2-k} \right) q_{ (n+2-k) }}^{q_{ (n+2-k) }} \Bigg(  \int\limits_{0}^{ q_{ (n+1-k) } } \dots \int\limits_{0}^{ q_{(2) } } dq_{ (1) }  \dots dq_{(n-k)}\Bigg) dq_{ (n+1-k) } \dots d q_{ (n-1) }  dq_{ (n) } $ \\[0.05 in ] 
$ = \prod\limits_{ j = n+1-k}^{ n-1} \Big( 1 - \left( \frac{j}{j+1} \right)^j \Big)$ \\[0.05 in] 
$ = \prod\limits_{ m=1}^{ k-1}\Big( 1 - \left( \frac{n-m}{n-m+1} \right)^{n-m} \Big) $ \\[0.05 in]
$ \leq (1 - e^{-1} )^{ k-1} $ \\[0.1 in] 
If we choose a sequence $ \{ x_n \} $ which is $o(n) $, and $ \sum\limits_{n=2}^{ \infty} (1-e^{-1})^{x_n - 1} < \infty $, then by first Borel Cantelli lemma, we can conclude that, $ P( j \leq x_n \text{ for all but finitely many n} ) = 1 $.  One such example might be $ x_n = n^{ 0.9} $. \\[0.05 in] 
So, $ j = o(n) $ almost surely. Note that, for any j with $ j = o (n) $,  $ \frac{ n+1-j}{ n (1- P_{ (j) } )} \geq 1 $ for all but finitely many n, with probability 1. Hence, 
$$ P( \hat{\pi}_{0,BH}   = 1 \text{ for all but finitely many n}) = 1  $$ 
This holds true not only for the global null hypothesis but also for the case where $\lim_{n \to \infty} \pi_{0} = 1$, which is encountered most often in real-life applications. To illustrate this, consider that $n_0$ of the $q_i$'s are derived from a uniform distribution $Unif [0,1]$, while the remaining originate from some alternate distribution. Let $q_{(1)}^{0} \leq \dots \leq q_{(n_0)}^0$ denote the order statistics among the sub-sample of null p-values. Let $j_0$ represent the index among the null p-values where $\frac{q_{(n_0 + 1 - k)}^0}{n_0 + 1 - k}$ first decreases. Then, in the combined sample, it is in the $( n_0 +1 - j_0 + x_1 )$-th place where $x_1 $ is the no. of non-null p values lying above that place. Clearly $$ (n_0 + 1 - j_0 + x_1) \leq n+1 - j $$ 
This implies, $$ j \leq j_0 + (n_1 - x_1) $$ 
By theorem \ref{pi_0=1}, $j_0 = o (n_0) = o (n) $ and since $ \lim_{ n \to \infty} \pi_0 = 1 $. this implies, $ (n_1 - x_1)  = o (n)$. Hence, $ j = o (n) $ and the conclusion of the theorem \ref{pi_0=1} remains valid. \\[0.05 in]
Now we ask the question, what is the behaviour of $\hat{\pi} _{0,BH} $ when $ 0 < \pi_0 < 1 $. 

\subsection{ Behaviour of $\hat{\pi}_{0,BH}$ for $0 < \pi_0 < 1 $} 
\begin{theorem}\label{general}
For any $0 < \pi_0 < 1 $, 
$$ \hat{\pi}_{0,BH} \stackrel{a.s.}{\to} 1 \text{ as } n \to \infty $$
\end{theorem} 
\textbf{Proof} : Let $h(.), H(.) $ denote the distribution of $q_i = 1 - P_i$ for $i = 1, \dots , n$. Then, the joint density of $q_{(1)} \leq q_{(2) } \leq \dots q_{ (n) } $  is given by $ h_n (x) = n ! ( \prod\limits_{i=1}^{n} h(x_i) )$. \\[0.05 in] 
Now, the proof proceeds in the similar manner as theorem \ref{pi_0=1}.  
Again, for $ k = 2, \dots , n $ \\[0.05 in] 
$ P( j > k ) $ \\[0.05 in]
$ = P( S_k \geq S_{k-1} \geq \dots \geq S_2 \geq S_1) $ \\[0.05 in]
$ = P \left( \frac{ q_{ (n+1-k) }}{ n+1-k } \geq  \frac{ q_{ (n+2-k) }}{ n+2-k } \geq  \frac{ q_{ (n) }}{ n } \right)$ \\[0.05 in] 
$ = 
n! \bigint\limits_{0}^{1} \bigint\limits_{ \left( \frac{n-1}{n} \right) q_{ (n) }}^{q_{ (n) }} \dots \bigint\limits_{ \left( \frac{n+1-k}{n+2-k} \right) q_{ (n+2-k) }}^{q_{ (n+2-k) }} \Bigg(  \int\limits_{0}^{ q_{ (n+1-k) } } \dots \int\limits_{0}^{ q_{(2) } } ( \prod\limits_{i=1}^{n} h(x_i) ) dq_{ (1) }  \dots dq_{(n-k)}\Bigg) dq_{ (n+1-k) } \dots d q_{ (n-1) }  dq_{ (n) } $ \\[0.05 in ] 
$ = \prod\limits_{ j = n+1-k}^{ n-1} \Big( 1 - \left( \frac{j}{j+1} \right)^j \Big)$ \text{  This is because } H(0) = 0 \text{ and } H(1) = 1  \\[0.05 in] 
$ = \prod\limits_{ m=1}^{ k-1}\Big( 1 - \left( \frac{n-m}{n-m+1} \right)^{n-m} \Big) $ \\[0.05 in]
It is very important to observe that, the quantity $j$ is actually a distribution-free quantity. So, we can again make a conclusion similar to theorem \ref{pi_0=1} that, $j = o(n) $ and consequently $$ \hat{\pi}_{0,BH} \stackrel{a.s.}{\to} 1 \text{ as } n \to \infty $$
However, the rate of convergence is very slow and depends on the value of $\pi_0$. We have considered testing $H_{0i} \:  : \: \mu_i = 0 $  vs  $ H_{1i} \: : \:  \mu_i > 0 $ in the context of Gaussian distribution with unit variance. Our approach involved drawing $n=10000$ observations from a blend of N(0,1) and N(2,1) distributions, and then deriving $\hat{\pi}_{0,BH}  $ based on the computed p-values. Repeating this process 10000 times allowed us to scrutinize the distribution of $\hat{\pi}_{0,BH} $. The ensuing table showcases the mean and minimum values of $\hat{\pi}_{0,BH} $ across these replications. Notably, for all instances where $ 0 < \pi_0 < 1$, the minimum estimators notably surpass the true $\pi_0$ value. It indicates that the adaptive BH procedure remains valid using the plug-in estimator. 
\begin{table}[h!]\caption{$\hat{\pi}_{0,BH}  $ for Gaussian mean test for $n = 10^4$}
\begin{center} 
\begin{tabular}{|l|l|l|l|l|}
\hline
$\pi_0$ & $E [\hat{\pi}_{0,BH}  ]$ & $ \text{ Median } \hat{\pi}_{0,BH}  $ & $\min  \hat{\pi}_{0,BH}  $ & Sd       \\ \hline
0.5     & 0.76               & 0.76                            & 0.69                 & 0.021    \\ \hline
0.75    & 0.91               & 0.91                            & 0.88                 & 0.011    \\ \hline
0.8     & 0.94               & 0.94                            & 0.91                 & 0.008    \\ \hline
0.9     & 0.98               & 0.98                            & 0.96                 & 0.004    \\ \hline
0.95    & 0.99               & 0.99                            & 0.98                 & 0.002    \\ \hline
1       & 1                  & 0.998                           & 1                    & 6.61E-05 \\ \hline
\end{tabular}

\end{center}
\end{table} 
\\
The table above illustrates that the value of $\hat{\pi}_{0,BH}$ significantly exceeds $\pi_0$ for $n = 10000$, confirming that the adaptive Benjamini-Hochberg (BH) procedure remains valid for practical applications. However, the primary result, $\hat{\pi}_{0,BH} \stackrel{a.s.}{\to} 1 \text{ as } n \to \infty$, is not evident even for $n = 10000$ due to the slow rate of convergence. Additionally, we conducted similar simulation exercises for a larger number of hypotheses ($n = 10^7, 10^9$). The subsequent table demonstrates that the BH estimator gradually converges to 1 as $n$ increases. The convergence rate is particularly slow when $\pi_0$ is significantly less than 1.
\begin{table}[ht] \caption{ $\hat{\pi}_{0,BH}  $ for Gaussian mean test for $n = 10^7$}
\begin{center}
\begin{tabular}{|l|l|l|l|l|}
\hline
$\pi_0$ & $ E [ \hat{\pi}_{0,BH} ] $ & $ \text{ Median} ( \hat{\pi}_{0,BH} ) $ & $ \text{ Min} ( \hat{\pi}_{0,BH} ) $ & Sd       \\ \hline
0.5     & 0.8428                     & 0.8417                                  & 0.8213                               & 0.0093   \\ \hline
0.75    & 0.9494                     & 0.9491                                  & 0.9392                               & 0.0037   \\ \hline
0.8     & 0.9648                     & 0.9646                                  & 0.9585                               & 0.0027   \\ \hline
0.9     & 0.9891                     & 0.9890                                  & 0.9866                               & 0.0011   \\ \hline
0.95    & 0.9968                     & 0.9968                                  & 0.9958                               & 0.0004   \\ \hline
1       & 1                          & 1                                       & 0.9999                           & 7.70E-08 \\ \hline
\end{tabular}
\end{center}
\end{table} 

\subsection{ Behaviour of $\hat{\pi}_{0,BH}$ under dependence } 
Our simulations have primarily focused on independent P-values. However, we have extended our analysis to explore the behavior of $\hat{\pi}_{0,BH}$ under various forms of dependence. Initially, we investigate the m-dependent scenario, where $X_i$ and $X_j$ are independent if $|i-j| > m$. In this setup, we generate observations $X_i$ from a set of independent normal variables $\xi_1, \xi_2, \dots $, following the formula $X_i = \frac{1}{\sqrt{m+1}} 
\sum\limits_{j = i}^{i+m} \xi_j $.

Within this framework, $n_1$ of the $X_i$'s are $N(\mu, 1)$ and $n_0$ of them are $N(0,1)$ random variables. For our simulations, we fix $n = 1000$ and explore six combinations of $\pi_0$, consistent with our previous investigations. Additionally, we vary the parameter $m$ with values of 2 and 5.
\begin{table}[ht]\caption{$\hat{\pi}_{0,BH}$ for Gaussian mean test for $m=2$ and $n=1000$}
\begin{center} 
\begin{tabular}{|l|l|l|l|l|}
\hline
$\pi_0$ & $ E [ \hat{\pi}_{0,BH} ] $ & $ \text{ Median} ( \hat{\pi}_{0,BH} ) $ & $ \text{ Min} ( \hat{\pi}_{0,BH} ) $ & Sd    \\ \hline
0.5     & 0.692                      & 0.688                                   & 0.600                                & 0.034 \\ \hline
0.75    & 0.903                      & 0.902                                   & 0.846                                & 0.019 \\ \hline
0.8     & 0.935                      & 0.934                                   & 0.888                                & 0.016 \\ \hline
0.9     & 0.983                      & 0.984                                   & 0.957                                & 0.007 \\ \hline
0.95    & 0.997                      & 0.998                                   & 0.983                                & 0.003 \\ \hline
1       & 1                          & 1                                       & 1                                    & 0     \\ \hline
\end{tabular}
\end{center} 
\end{table} 

\begin{table}[ht]
\caption{$\hat{\pi}_{0,BH}$ for Gaussian mean test for $m=5$ and $n=1000$}
\begin{center} 
\begin{tabular}{|l|l|l|l|l|}
\hline
$\pi_0$ & $ E [ \hat{\pi}_{0,BH} ] $ & $ \text{ Median} ( \hat{\pi}_{0,BH} ) $ & $ \text{ Min} ( \hat{\pi}_{0,BH} ) $ & Sd    \\ \hline
0.5     & 0.621                      & 0.614                                   & 0.574                                & 0.038 \\ \hline
0.75    & 0.896                      & 0.882                                   & 0.826                                & 0.044 \\ \hline
0.8     & 0.943                      & 0.934                                   & 0.875                                & 0.035 \\ \hline
0.9     & 0.996                      & 1.000                                   & 0.959                                & 0.007 \\ \hline
0.95    & 1.000                      & 1.000                                   & 0.989                                & 0.001 \\ \hline
1       & 1                          & 1                                       & 1                                    & 0     \\ \hline
\end{tabular}
\end{center} 
\end{table}
It is interesting to observe that, in all cases, the value of $\hat{\pi}_{0,BH}$ lies well above the actual value. We have further simulated the same for $n=10000$ and the observations are presented below. 
\begin{table}[h!]
\caption{$\hat{\pi}_{0,BH}$ for Gaussian mean test for $m=2$ and $n=10000$}
\begin{center} 
\begin{tabular}{|l|l|l|l|l|}
\hline
$\pi_0$ & $ E [ \hat{\pi}_{0,BH} ] $ & $ \text{ Median} ( \hat{\pi}_{0,BH} ) $ & $ \text{ Min} ( \hat{\pi}_{0,BH} ) $ & Sd    \\ \hline
0.5     & 0.746                      & 0.742                                   & 0.681                                & 0.027 \\ \hline
0.75    & 0.931                      & 0.930                                   & 0.896                                & 0.013 \\ \hline
0.8     & 0.956                      & 0.955                                   & 0.926                                & 0.010 \\ \hline
0.9     & 0.991                      & 0.991                                   & 0.978                                & 0.004 \\ \hline
0.95    & 0.999                      & 0.999                                   & 0.993                                & 0.001 \\ \hline
1       & 1                          & 1                                       & 1                                    & 0     \\ \hline
\end{tabular}
\end{center} 
\end{table}  
\begin{table}[h!]
\centering
\caption{$\hat{\pi}_{0,BH}$ for Gaussian mean test for $m=5$ and $n=10000$}
\begin{tabular}{|l|l|l|l|l|}
\hline
$\pi_0$ & $ E [ \hat{\pi}_{0,BH} ] $ & $ \text{ Median} ( \hat{\pi}_{0,BH} ) $ & $ \text{ Min} ( \hat{\pi}_{0,BH} ) $ & Sd       \\ \hline
0.5     & 0.878                      & 1                                  & 0.608                                & 0.163    \\ \hline
0.75    & 0.998                      & 1                                       & 0.880                                & 0.011    \\ \hline
0.8     & 0.999                      & 1                                       & 0.935                                & 0.002    \\ \hline
0.9     & 1                          & 1                                       & 0.999                                & 2.04E-05 \\ \hline
0.95    & 1                          & 1                                       & 1                                    & 0        \\ \hline
1       & 1                          & 1                                       & 1                                    & 0        \\ \hline
\end{tabular}
\end{table} 
\newpage

As before, the values of $\hat{\pi}_{0,BH}$ slowly increases as we increase $n$. So, we might expect similar result for m-dependent setup also.\\[0.05 in]
Now, we explore block dependent structure. First, we have considered the problem of Gaussian mean test with 10 blocks of equicorrelated normal variables with correlations $0.1, \dots , 0.9, 0.95$ and each block consists of $100$ observations. The observations are presented in the following table. 
\begin{table}[h!]
\centering
\caption{$\hat{\pi}_{0,BH}$ for Gaussian mean test under block dependence for $n=1000$}
\centering
\begin{tabular}{|l|l|l|l|l|}
\hline
$\pi_0$ & $ E [ \hat{\pi}_{0,BH} ] $ & $ \text{ Median} ( \hat{\pi}_{0,BH} ) $ & $ \text{ Min} ( \hat{\pi}_{0,BH} ) $ & Sd    \\ \hline
0.5     & 0.713                      & 0.712                                   & 0.349                                & 0.082 \\ \hline
0.75    & 0.886                      & 0.891                                   & 0.547                                & 0.051 \\ \hline
0.8     & 0.915                      & 0.922                                   & 0.635                                & 0.044 \\ \hline
0.9     & 0.965                      & 0.972                                   & 0.690                                & 0.030 \\ \hline
0.95    & 0.985                      & 0.991                                   & 0.745                                & 0.023 \\ \hline
1       & 0.996                      & 1                                       & 0.765                                & 0.018 \\ \hline
\end{tabular}
\end{table}
Although the mean and median of the $10000$ replications lie well above the original value of $\pi_0$, the minimum observed from the $10000$ replications is not above $\pi_0$ anymore. This is possibly due to the increased level of dependence among the hypotheses. However, this gap appears to narrow progressively and the estimator increases slowly as we continually increase $n$ from $1000$ to $10000$.
\begin{table}[h!]
\centering
\caption{$\hat{\pi}_{0,BH}$ for Gaussian mean test under block dependence for $n=10000$}
\centering
\begin{tabular}{|l|l|l|l|l|}
\hline
$\pi_0$ & $ E [ \hat{\pi}_{0,BH} ] $ & $ \text{ Median} ( \hat{\pi}_{0,BH} ) $ & $ \text{ Min} ( \hat{\pi}_{0,BH} ) $ & Sd    \\ \hline
0.5     & 0.761                      & 0.761                                   & 0.480                                & 0.076 \\ \hline
0.75    & 0.911                      & 0.915                                   & 0.624                                & 0.042 \\ \hline
0.8     & 0.935                      & 0.940                                   & 0.688                                & 0.035 \\ \hline
0.9     & 0.975                      & 0.979                                   & 0.717                                & 0.022 \\ \hline
0.95    & 0.990                      & 0.993                                   & 0.773                                & 0.016 \\ \hline
1       & 0.998                      & 1                                       & 0.762                                & 0.011 \\ \hline
\end{tabular}
\end{table}

\section{Estimator proposed in \cite{storey2002direct}} 
Since the largest p-values are most likely to be uniformly distributed, \cite{storey2002direct} suggested that a conservative estimator of $\pi_0$ is 
\begin{equation}\label{storey}
 \hat{\pi}_{0}  ( \lambda ) = \frac{ \# \{ P_i > \lambda \} }{n (1- \lambda) } = \frac{ W ( \lambda )}{ n (1 - \lambda ) }
\end{equation} 
for some well chosen $\lambda$. Since $ Pr( P_i > \lambda ) =  \pi_0 (1- \lambda )  + \pi_1 (1 - F_1 ( \lambda ))  $, by strong law of large numbers, one can easily deduce that for any $0 < \lambda < 1 $ and for independent p values, 
$$
\hat{\pi}_{0}  ( \lambda) \stackrel{a.s.}{\to}  \pi_0 + \pi_1 \left( \frac{1 - F_1 ( \lambda) }{1 - \lambda } \right) \text{  as } n \to \infty $$
So, the estimator is asymptotically a conservative estimator of $\pi_0$. \cite{storey2002direct} proposed two methods for controlling pFDR and FDR and showcased the gain in power and accuracy using the estimator $\hat{\pi}_{0}  ( \lambda) $. While the validity of the results were initially established under the assumption of independent p-values, it can be shown that these findings remain applicable even under conditions of weak dependence. For p values with empirical cdf and cdf $F_n$ and $F$ respectively, if $|| F_n - F||_{ \infty} \stackrel{a.s.}{\to} 0 $ as $n \to \infty$, then all theorems of \cite{storey2002direct} remains valid. \\ 
Under various forms of weak dependence among p-values, the condition $||F_n - F||_{ \infty} \stackrel{a.s.}{\to}  0 $ persists. By invoking the Ergodic theorem as presented in \cite{billingsley2017probability}, it becomes apparent that if the right shift transformation of the sequence of p-values is measure-preserving (i.e., $(P_1, P_2, P_3, \dots ) \stackrel{d}{=} (P_2 , P_3 , \dots ) $), then the empirical cumulative distribution function (ecdf) $ F_n (x) $ converges pointwise almost surely to some limit as $ n \to \infty $. The condition of measure preservation in the right shift transformation is synonymous with stating that the sequence of p-values is strongly stationary. The stationarity of the p-values guarantees the existence of an almost sure limit, which may or may not be degenerate.\\[0.05 in]
To ensure that the pointwise limit equals to the cdf $F(x) $, it suffices to assume that $ F_n (x) \stackrel{P}{\to} F(x) \: \forall \: x $. Convergence in probability implies that a subsequence of the original sequence is almost surely convergent. Given the original sequence's almost sure convergence, all subsequences should converge to the same limit, thereby leading to $F_n (x) $ converging almost surely to $F(x) $ for all $x$. Consequently, according to the Glivenko-Cantelli lemma, it can be asserted that $|| F_n - F||_{ \infty} \stackrel{a.s.}{\to} 0 $ as $n \to \infty$. These discussions are formally encapsulated into the theorem presented below.
 \begin{tcolorbox}
 \begin{theorem}\label{EMP_CDF}
Let $P_1, P_2, \dots , P_n $ be a sequence of p values in a multiple testing problem with empirical cdf and continuous cdf $F_n , F $ respectively. We assume the following about the joint distribution of the p values.
\begin{enumerate}[label=\Roman*]
    \item) $ \; \; P_1, P_2 , ...$ forms a strongly stationary sequence. 
    \item) $ \; \; F_n (x) \stackrel{P}{\to} F(x) \: \forall \: x $ as $ n \to \infty $.  
\end{enumerate} 
Then, $|| F_n - F||_{ \infty} \stackrel{a.s.}{\to} 0 $ as $n \to \infty$
\end{theorem}
\end{tcolorbox}  
The assumption (II) outlined in theorem \ref{EMP_CDF} holds over a wide array of dependence structures commonly encountered in practical contexts. We highlight a few such cases where the assumptions are met. 
\begin{itemize}
    \item m-dependent and stationary p values satisfy assumption (I) and (II) 
    \item Hidden Markov model similar to \cite{sun2009large}. In particular. if the sequence of P values form a stationary, Ergodic Markov chain, then the assumptions (I) and (II) are met. 
    \item All Ergodic, stationary sequence of p values satisfy assumption (I) and (II). This is because, ergodicity ensures that the limit of $F_n(.) $ is indeed the distribution function. In practice, we often assume a sequence of random variables to be asymptotically independent. One such example of asymptotic independence is mixing sequence (\cite{billingsley2017probability}). Since all mixing sequences are ergodic (under the right shift operation), mixing sequences also satisfy assumption (I) and (II). 
    \item The assumption of asymptotic independence can be somewhat relaxed as follows. 
    $$ \forall \: x, \: \: \lim_{ n \to \infty} \frac{1}{n^2} 
   \underset{i \neq j }{\sum \sum}  (Pr( P_i \leq x , P_j \leq x ) - (F(x))^2) = 0 $$ 
    By Chebyshev's inequality, it can be easily argued that, under the above mentioned assumption, $ F_n (x) \stackrel{P}{\to} F(x) \: \forall \: x $ as $ n \to \infty $. In fact, it constitutes a necessary and sufficient criterion for assumption (II) to hold. 
\end{itemize}  
\subsection{Simulation results}
The estimator in \cite{storey2003positive} involves a tuning parameter, $\lambda$. Choosing $\lambda$ requires balancing bias and variance for the estimator $\hat{\pi}_0(\lambda)$. According to \cite{storey2003positive}, for well-behaved p-values, the bias tends to decrease as $\lambda$ increases, reaching its minimum as $\lambda$ approaches 1. Consequently, \cite{storey2003positive} proposed a method for selecting the optimal $\lambda$. 
\begin{itemize}
    \item Fix a range of $\lambda$'s by $ \Lambda $ (e.g. $ \Lambda = \{ 0 , 0.01, 0.02 , \dots , 0.95 \} $) 
    \item By \ref{storey}, compute $ \hat{\pi}_0 ( \lambda ) $ for each $ \lambda \in \Lambda$. 
    \item Fit a cubic spline ($\hat{f}$) on $ \hat{\pi}_0 ( \lambda ) $ and obtain the final estimate as 
    $$ \hat{ \pi}_0 = \min \{ \hat{f} (1) , 1 \}  $$
\end{itemize}
In addition to the smoothing-based technique introduced by \cite{storey2003positive}, a bootstrap-based method for selecting the optimal value of \( \lambda \) was proposed in \cite{storey2004strong}. This method builds on the earlier work of \cite{storey2002direct}. The proposed automatic choice of \( \lambda \) aims to estimate the value that minimizes the mean squared error (MSE), balancing bias and variance. Specifically, it seeks to minimize \( E \left[ (\hat{\pi}_0(\lambda) - \pi_0)^2 \right] \). The procedure for this method is summarized below.
\begin{itemize}
    \item For each $\lambda \in \Lambda$, compute $\hat{\pi}_0 
 ( \lambda ) $ as \ref{storey}. 
 \item Generate $B$ bootstrap samples from the p values. For each $ b = 1 , \dots , B$, find the estimators $ \{\hat{\pi}^{ *b} _0 ( \lambda )  \}_{ \lambda \in \Lambda } $ based on the $b$-th sample.  
 \item Since $ E [ \hat{ \pi}_0  ( \lambda ) ] \geq \\pi_0 $ for all $ \lambda \in [0, 1 ) $, a plug-in estimator of $ \pi_0$ can be taken as 
 $$ \hat{\pi}_0 = \min_{ \lambda^{'} \in \Lambda} \{ \hat{ \pi}_0 ( \lambda^{'} ) \} $$
 \item for each $\lambda \in \Lambda $, estimate its respective MSE as 
 $$ \hat{MSE} ( \lambda ) = \frac{1}{B} \sum\limits_{b=1}^{B} 
 \left[  \hat{\pi}^{ *b} _0 ( \lambda )  - \min_{ \lambda^{'} \in \Lambda} \{ \hat{ \pi}_0 ( \lambda^{'} ) \} \right]^2
 $$
 \item Set the final estimator as 
 $$ \hat{\pi }_0 = \min \left\{ 1 , \text{argmin}_{ \lambda \in \Lambda } \hat{MSE} ( \lambda ) \right\} $$
\end{itemize}
Based on these two methods, we examined the performance of Storey's estimator under various dependence structures. Initially, we focused on the autoregressive model of order 1 (AR(1)), where observations are generated according to the mechanism described below.
\begin{equation}\label{AR1}
    X_i = a X_{i-1} + \xi_i \; \:  \text{ for } i > 1 \text{ and } X_1 = \xi_1  
    \end{equation}
$\xi_i$'s are standard normal random variables. Theorem 19.1 of \cite{billingsley2013convergence} ensures that the AR(1) process meets the necessary criteria for the validity of Storey's estimator. In our analysis, we set \( n = 1000 \) and considered values of \( a \) as \( 0.1, \: 0.2, \: 0.5, \: 0.75 \). We also examined 10 different combinations of \( \pi_0 \), specifically \( \{0.5, \: 0.6, \: 0.65, \: 0.7, \: \dots, \: 0.95, \: 0.98\} \). Observations were generated from a mixture of \( N(0,1) \) and \( N(2,1) \) random variables, and p-values were computed for a one-sided test of the mean. \\[0.04 in]
For values of \( a = 0.1 \) and \( 0.2 \), both the bootstrap and spline-based estimators consistently yielded results above the true values of \( \pi_0 \). There was minimal difference between these two methods in terms of selecting the optimal \( \lambda \). 

\newpage

\begin{center}
   \begin{table}[h!]
\centering
\caption{$n = 1000$ and $a = 0.1, \: 0.2 $}
\label{a0102}
\begin{tabular}{|l|ll|ll|}
\hline
\textbf{}        & \multicolumn{2}{c|}{\textbf{$ a = 0.1$}}                                                      & \multicolumn{2}{c|}{\textbf{$ a = 0.2$}}                                                      \\ \hline
$\mathbf{\pi_0} $ & \multicolumn{1}{l|}{\textbf{$ \hat{\pi}_0$ (Bootstrap)}} & \textbf{$ \hat{\pi}_0$ (Smoother)} & \multicolumn{1}{l|}{\textbf{$ \hat{\pi}_0$ (Bootstrap)}} & \textbf{$ \hat{\pi}_0$ (Smoother)} \\ \hline
0.5              & \multicolumn{1}{l|}{0.501}                               & 0.502                              & \multicolumn{1}{l|}{0.522}                               & 0.522                              \\ \hline
0.6              & \multicolumn{1}{l|}{0.604}                               & 0.605                              & \multicolumn{1}{l|}{0.628}                               & 0.627                              \\ \hline
0.65             & \multicolumn{1}{l|}{0.655}                               & 0.655                              & \multicolumn{1}{l|}{0.681}                               & 0.680                              \\ \hline
0.7              & \multicolumn{1}{l|}{0.707}                               & 0.707                              & \multicolumn{1}{l|}{0.734}                               & 0.735                              \\ \hline
0.75             & \multicolumn{1}{l|}{0.756}                               & 0.757                              & \multicolumn{1}{l|}{0.785}                               & 0.786                              \\ \hline
0.8              & \multicolumn{1}{l|}{0.809}                               & 0.808                              & \multicolumn{1}{l|}{0.838}                               & 0.839                              \\ \hline
0.85             & \multicolumn{1}{l|}{0.858}                               & 0.858                              & \multicolumn{1}{l|}{0.892}                               & 0.892                              \\ \hline
0.9              & \multicolumn{1}{l|}{0.911}                               & 0.908                              & \multicolumn{1}{l|}{0.945}                               & 0.944                              \\ \hline
0.95             & \multicolumn{1}{l|}{0.959}                               & 0.959                              & \multicolumn{1}{l|}{0.987}                               & 0.985                              \\ \hline
0.98             & \multicolumn{1}{l|}{0.984}                               & 0.984                              & \multicolumn{1}{l|}{0.997}                               & 0.997                              \\ \hline
\end{tabular}
\end{table}
\end{center}

It is interesting to note that as the degree of dependence increases, Storey's estimator tends to become more conservative. Specifically, as we move from \( a = 0.1 \) and \( a = 0.2 \) to \( a = 0.5 \) and \( a = 0.75 \), the estimated values approach 1 for all values of \( \pi_0 \).

\begin{center}
    \begin{table}[ht]
\centering
\caption{$n = 1000$ and $a = 0.5, \: 0.75$}
\label{a05075}
\begin{tabular}{|l|ll|ll|}
\hline
\textbf{}        & \multicolumn{2}{c|}{\textbf{$ a = 0.5$}}                                                      & \multicolumn{2}{c|}{\textbf{$ a = 0.75$}}                                                     \\ \hline
\textbf{$\pi_0$} & \multicolumn{1}{l|}{\textbf{$ \hat{\pi}_0$ (Bootstrap)}} & \textbf{$ \hat{\pi}_0$ (Smoother)} & \multicolumn{1}{l|}{\textbf{$ \hat{\pi}_0$ (Bootstrap)}} & \textbf{$ \hat{\pi}_0$ (Smoother)} \\ \hline
0.5              & \multicolumn{1}{l|}{0.692}                               & 0.691                              & \multicolumn{1}{l|}{0.977}                               & 0.977                              \\ \hline
0.6              & \multicolumn{1}{l|}{0.827}                               & 0.828                              & \multicolumn{1}{l|}{0.998}                               & 0.996                              \\ \hline
0.65             & \multicolumn{1}{l|}{0.893}                               & 0.895                              & \multicolumn{1}{l|}{0.998}                               & 0.997                              \\ \hline
0.7              & \multicolumn{1}{l|}{0.962}                               & 0.961                              & \multicolumn{1}{l|}{0.999}                               & 1                                  \\ \hline
0.75             & \multicolumn{1}{l|}{0.995}                               & 0.995                              & \multicolumn{1}{l|}{1}                                   & 1                                  \\ \hline
0.8              & \multicolumn{1}{l|}{1}                                   & 1                                  & \multicolumn{1}{l|}{1}                                   & 1                                  \\ \hline
0.85             & \multicolumn{1}{l|}{1}                                   & 1                                  & \multicolumn{1}{l|}{1}                                   & 1                                  \\ \hline
0.9              & \multicolumn{1}{l|}{1}                                   & 1                                  & \multicolumn{1}{l|}{1}                                   & 1                                  \\ \hline
0.95             & \multicolumn{1}{l|}{1}                                   & 1                                  & \multicolumn{1}{l|}{1}                                   & 1                                  \\ \hline
0.98             & \multicolumn{1}{l|}{1}                                   & 1                                  & \multicolumn{1}{l|}{1}                                   & 1                                  \\ \hline
\end{tabular}
\end{table}
\end{center}
In the AR(1) model, Storey's estimator exhibits increased bias as the correlation strengthens. Additionally, we evaluated the estimator for a stationary \( m \)-dependent sequence of p-values. In this case, we generated p-values from a mixture of normal distributions using the following approach: \\[0.05 in] 
Let $ \xi_1, \xi_2, \dots \stackrel{i.i.d.}{\sim} N(0,1) $ and $ X_i  \propto \sum\limits_{ j = i}^{ i+m } \left( (m+1) - (j-i) \right) \xi_j $. $X_i$'s are standardized so that, they have unit variance. In this case, $\{ X_i \}_{ i \geq 1 } $ forms a m-dependent sequence of random variables with the nearer observations getting higher weight. Next, we add 2 to \( n (1 - \pi_0) \) of the \( X_i \) values so that they form a mixture of \( N(0,1) \) and \( N(2,1) \) random variables. When computing Storey's estimator based on these p-values, the estimates tend to be very close to 1, indicating a higher level of dependence. 
\begin{center}
\begin{table}[ht]
\centering
\caption{$n = 1000$ and $m=1,2$}
\label{m12}
\begin{tabular}{|l|ll|ll|}
\hline
\textbf{}       & \multicolumn{2}{c|}{\textbf{m=1}}                                                               & \multicolumn{2}{c|}{\textbf{m=2}}                                                               \\ \hline
\textbf{$pi_0$} & \multicolumn{1}{l|}{\textbf{$ \hat{\pi}_0 $ (Bootstrap)}} & \textbf{$ \hat{\pi}_0 $ (Smoother)} & \multicolumn{1}{l|}{\textbf{$ \hat{\pi}_0 $ (Bootstrap)}} & \textbf{$ \hat{\pi}_0 $ (Smoother)} \\ \hline
0.5             & \multicolumn{1}{l|}{0.915}                                & 0.915                               & \multicolumn{1}{l|}{0.999}                                & 0.999                               \\ \hline
0.6             & \multicolumn{1}{l|}{0.993}                                & 0.992                               & \multicolumn{1}{l|}{1}                                    & 1                                   \\ \hline
0.65            & \multicolumn{1}{l|}{0.999}                                & 0.999                               & \multicolumn{1}{l|}{1}                                    & 1                                   \\ \hline
0.7             & \multicolumn{1}{l|}{0.999}                                & 0.999                               & \multicolumn{1}{l|}{1}                                    & 1                                   \\ \hline
0.75            & \multicolumn{1}{l|}{1}                                    & 1                                   & \multicolumn{1}{l|}{1}                                    & 1                                   \\ \hline
0.8             & \multicolumn{1}{l|}{1}                                    & 1                                   & \multicolumn{1}{l|}{1}                                    & 1                                   \\ \hline
0.85            & \multicolumn{1}{l|}{1}                                    & 1                                   & \multicolumn{1}{l|}{1}                                    & 1                                   \\ \hline
0.9             & \multicolumn{1}{l|}{1}                                    & 1                                   & \multicolumn{1}{l|}{1}                                    & 1                                   \\ \hline
0.95            & \multicolumn{1}{l|}{1}                                    & 1                                   & \multicolumn{1}{l|}{1}                                    & 1                                   \\ \hline
0.98            & \multicolumn{1}{l|}{1}                                    & 1                                   & \multicolumn{1}{l|}{1}                                    & 1                                   \\ \hline
\end{tabular}
\end{table}
\end{center} 
And if we move further to $m=5$, then the estimator is identically 1 for all choices of $\pi_0$. The simulation results vary slightly in the $m$-dependent case, where observations are generated as 
$
X_i = \frac{1}{\sqrt{m+1}} \sum\limits_{j=i}^{i+m} \xi_j.
$
Unlike scenarios with decreasing weights, this setup uses equal weights for the $m$ neighboring observations. The heightened dependence in this case results in slower convergence. Consequently, Storey's estimator is much closer to the true value of $\pi_0$ compared to the previous case, where estimated values were near 1. Although Storey's estimator is typically biased and tends to overestimate $\pi_0$, this effect is less pronounced here due to the increased dependence among observations. 
\newpage
\begin{center}
 \begin{table}[ht]
\centering
\caption{$n = 1000$ and $m=1,2,5$ for the equally weighted case}
\label{m-dependent_Storey}
\resizebox{\columnwidth}{!}{%
\begin{tabular}{|l|ll|ll|ll|}
\hline
\textbf{} & \multicolumn{2}{l|}{\textbf{m=1}}  & \multicolumn{2}{l|}{\textbf{m=2}}  & \multicolumn{2}{l|}{\textbf{m=5}}  \\ \hline
\textbf{$\pi_0$} &
  \multicolumn{1}{l|}{\textbf{$\hat{\pi}_0$ (Bootstrap)}} &
  \textbf{$\hat{\pi}_0$ (Smoother)} &
  \multicolumn{1}{l|}{\textbf{$\hat{\pi}_0$ (Bootstrap)}} &
  \textbf{$\hat{\pi}_0$ (Smoother)} &
  \multicolumn{1}{l|}{\textbf{$\hat{\pi}_0$ (Bootstrap)}} &
  \textbf{$\hat{\pi}_0$ (Smoother)} \\ \hline
0.5       & \multicolumn{1}{l|}{0.494} & 0.497 & \multicolumn{1}{l|}{0.496} & 0.497 & \multicolumn{1}{l|}{0.501} & 0.494 \\ \hline
0.6       & \multicolumn{1}{l|}{0.598} & 0.599 & \multicolumn{1}{l|}{0.600} & 0.597 & \multicolumn{1}{l|}{0.591} & 0.597 \\ \hline
0.65      & \multicolumn{1}{l|}{0.642} & 0.650 & \multicolumn{1}{l|}{0.647} & 0.647 & \multicolumn{1}{l|}{0.649} & 0.644 \\ \hline
0.7       & \multicolumn{1}{l|}{0.694} & 0.700 & \multicolumn{1}{l|}{0.699} & 0.704 & \multicolumn{1}{l|}{0.687} & 0.701 \\ \hline
0.75      & \multicolumn{1}{l|}{0.755} & 0.749 & \multicolumn{1}{l|}{0.742} & 0.744 & \multicolumn{1}{l|}{0.747} & 0.753 \\ \hline
0.8       & \multicolumn{1}{l|}{0.800} & 0.791 & \multicolumn{1}{l|}{0.797} & 0.797 & \multicolumn{1}{l|}{0.795} & 0.794 \\ \hline
0.85      & \multicolumn{1}{l|}{0.844} & 0.846 & \multicolumn{1}{l|}{0.839} & 0.841 & \multicolumn{1}{l|}{0.820} & 0.837 \\ \hline
0.9       & \multicolumn{1}{l|}{0.893} & 0.887 & \multicolumn{1}{l|}{0.884} & 0.880 & \multicolumn{1}{l|}{0.872} & 0.866 \\ \hline
0.95      & \multicolumn{1}{l|}{0.927} & 0.922 & \multicolumn{1}{l|}{0.919} & 0.918 & \multicolumn{1}{l|}{0.906} & 0.897 \\ \hline
0.98      & \multicolumn{1}{l|}{0.945} & 0.946 & \multicolumn{1}{l|}{0.936} & 0.934 & \multicolumn{1}{l|}{0.911} & 0.912 \\ \hline
\end{tabular}
}
\end{table}
\end{center}  

Due to the slight underestimation of $\pi_0$, the validity of the estimator seems to be questionable for $n=1000$. However, as $n$ increases, the difference appears to diminish. The bootstrap-based estimator exhibits a particularly interesting property in this context: it seems to converge to the true value of $\pi_0$, as evidenced by the simulation results for $n=10^6$. For all values of $m$ (i.e., $m=1,2,5$), the estimated values, rounded to three decimal places, show identical figures, as is apparent from Table~\ref{m-dependent_Storey_2}.

\begin{center}
    \begin{table}[ht]
\centering
\caption{$n = 10^6$ and $m=1,2,5$ for the equally weighted case}
\label{m-dependent_Storey_2}
\begin{tabular}{|l|l|l|l|}
\hline
     & \textbf{m=1} & \textbf{m=2} & \textbf{m=5} \\ \hline
\textbf{$\pi_0$} & \textbf{$\hat{\pi}_0$ (Bootstrap)} & \textbf{$\hat{\pi}_0$ (Bootstrap)} & \textbf{$\hat{\pi}_0$ (Bootstrap)} \\ \hline
0.5  & 0.496        & 0.496        & 0.496        \\ \hline
0.6  & 0.597        & 0.597        & 0.597        \\ \hline
0.65 & 0.647        & 0.647        & 0.647        \\ \hline
0.7  & 0.698        & 0.698        & 0.698        \\ \hline
0.75 & 0.748        & 0.748        & 0.748        \\ \hline
0.8  & 0.798        & 0.798        & 0.798        \\ \hline
0.85 & 0.849        & 0.849        & 0.849        \\ \hline
0.9  & 0.899        & 0.899        & 0.899        \\ \hline
0.95 & 0.950        & 0.950        & 0.950        \\ \hline
0.98 & 0.980        & 0.980        & 0.980        \\ \hline
\end{tabular}
\end{table}
\end{center}

A block-dependent sequence does not meet the stationarity requirements necessary for establishing the validity of the estimator via the Strong Law of Large Numbers (SLLN). We have examined a scenario with \( n = 1000 \), where the data is divided into 10 equal-sized blocks of equicorrelated standard normal random variables with correlations $ 0.1 , \: 0.2 , \: \dots , 0.9 , \: 0.95 $. \\[0.04 in]
Our analysis indicates that the estimates of $\pi_0$ are quite accurate for smaller values of $\pi_0$. However, as $\pi_0$ increases, particularly for values above $0.9$, the discrepancy between the estimated values and the true value of $\pi_0$ becomes more pronounced.
\begin{center}
\begin{table}[ht]
\centering
\caption{Block dependent variables with $n=1000$ and 10 blocks }
\label{block_m100}
\begin{tabular}{|l|l|l|l|l|l|l|l|l|l|l|}
\hline
\textbf{$\pi_0$}                   & 0.5   & 0.6   & 0.65  & 0.7   & 0.75  & 0.8   & 0.85  & 0.9   & 0.95  & 0.98  \\ \hline
\textbf{$\hat{\pi}_0$ (Bootstrap)} & 0.497 & 0.596 & 0.649 & 0.698 & 0.745 & 0.780 & 0.827 & 0.856 & 0.884 & 0.899 \\ \hline
\textbf{$\hat{\pi}_0$ (Smoother)}  & 0.497 & 0.602 & 0.644 & 0.693 & 0.737 & 0.785 & 0.827 & 0.871 & 0.886 & 0.909 \\ \hline
\end{tabular}
\end{table}
\end{center} 
Even with a slight reduction in dependence by considering $200$ blocks of equicorrelated variables, each of size $5$, the overall picture remains largely unchanged. In this scenario, the discrepancy between the estimated and true values of $\pi_0$ still becomes more pronounced for higher values of $\pi_0$.
\begin{center}
    \begin{table}[ht]
\centering
\caption{Block dependent variables with $n=1000$ and 200 blocks }
\label{block_m5}
\begin{tabular}{|l|l|l|l|l|l|l|l|l|l|l|}
\hline
\textbf{$\pi_0$}                   & 0.5   & 0.6   & 0.65  & 0.7   & 0.75  & 0.8   & 0.85  & 0.9   & 0.95  & 0.98  \\ \hline
\textbf{$\hat{\pi}_0$ (Bootstrap)} & 0.492 & 0.590 & 0.649 & 0.692 & 0.746 & 0.782 & 0.822 & 0.863 & 0.891 & 0.909 \\ \hline
\textbf{$\hat{\pi}_0$ (Smoother)}  & 0.501 & 0.590 & 0.645 & 0.698 & 0.742 & 0.786 & 0.829 & 0.866 & 0.886 & 0.911 \\ \hline
\end{tabular}
\end{table}
\end{center}
\section{Estimator proposed in \cite{patra2016estimation} }
If $\pi_1$ were known, for an independent and identically distributed sample $X_1, X_2 , \dots , X_n$ with distribution function $F$, a naive estimator of $F_1$ would be 
$$ \hat{F}_{1,n} = \frac{ F_n - \pi_0 F_0}{\pi_1} $$ 
where $F_n$ is the empirical CDF of the observed sample. Since $\hat{F}_{1,n} $ need not be a valid cdf, this naive estimator can be improved further as 
$$ \hat{F}_{2,n} = \text{\stackunder{argmin}{$ W \in \mathcal{F}$} } \: \: \frac{1}{n} \sum\limits_{i=1}^{n} \{ W(X_i) - \hat{F}_{1,n} \}^2 $$ 
where $\mathcal{F}$ is the set of all cdfs. The minimizer can be easily found by the pool adjacent violators algorithm mentioned in \cite{robertson1988order}. Using these two estimators $\hat{F}_{1,n}$ and $\hat{F}_{2,n}$, \cite{patra2016estimation} suggested the following estimator of $\pi_0$. 
\begin{equation}\label{est} 
\centering
\hat{\pi}_{1}^{c_n} = \inf \left\{  \gamma \in (0,1] \: : \: \gamma d_n ( \hat{F}_{1,n}, \hat{F}_{2,n} ) \leq \frac{c_n}{\sqrt{n}} \right\}
\end{equation}
where $c_n$ is a sequence of constants satsifying $c_n = o( \sqrt{n})  $  and $d_n (.,.)$ is same as the $ \mathcal{L}_2 ( F_n) $ distance, i.e. if $f,g \: : \: \mathbb{R} \to \mathbb{R} $ are two functions, then 
$$ d_{n}^2 (f,g) = \int ( f(x) - g(x))^2 d F_n (x) = \frac{1}{n} \sum\limits_{i =1}^{n} ( f(X_{(i) } )  - g( X_{ (i) } ) )^2 $$
The idea of this estimator originated from the following observations. 
\begin{itemize}
\item If $\gamma = 1 $, then clearly $ \hat{F}_{1,n} = \hat{F}_{2,n} = F_n $ and they are much different for smaller values of $\gamma$. 
\item For any $ \gamma \in (0,1]$, it is easy to see that, 
$$ d_n ( F_n , \gamma \hat{F}_{2,n} + (1- \gamma) F_0) = \gamma d_n ( \hat{F}_{1,n}, \hat{F}_{2,n} ) $$
And 
$$ \lim_{\gamma \to 0 }\gamma d_n ( \hat{F}_{1,n}, \hat{F}_{2,n} ) = d_n ( F_n , F_0)  $$

\item  $\gamma d_n ( \hat{F}_{1,n}, \hat{F}_{2,n} ) $ changes its behaviour at  \underbar{$\pi_1$}. 
\begin{theorem}\label{th_cons}
   For    \underbar{$\pi_1$} $\leq \gamma \leq 1$, $ \gamma d_n ( \hat{F}_{1,n}, \hat{F}_{2,n} ) \leq d_n ( F , F_n) $. Then, if $|| F_n - F ||_{ \infty} \stackrel{a.s.}{\to}  0 $ as $ n \to \infty$ then, 
   \begin{equation*}
\gamma d_n ( \hat{F}_{1,n}, \hat{F}_{2,n} ) \stackrel{a.s.}{\to}  
\begin{cases} 
0 \qquad \text{  if } \gamma \geq \text{\underbar{$\pi_1$}} \\
> 0, \: \:  \text{  if } \gamma < \text{\underbar{$\pi_1$}}
\end{cases}
\end{equation*}
\end{theorem} 
\end{itemize} 
The proof of Theorem \ref{th_cons} can be easily done using the idea similar to \cite{patra2016estimation}. The proof is given in appendix. Theorem \ref{th_cons} demonstrates that the function $\gamma d_n ( \hat{F}_{1,n}, \hat{F}_{2,n} )$ changes behaviour at $\underline{\pi_1}$. For sufficiently large $n$, the function has a steep decline to the left of $\underline{\pi_1}$ and a more gradual decrease near zero to the right of $\underline{\pi_1}$. Thus, $\gamma d_n ( \hat{F}_{1,n}, \hat{F}_{2,n} )$ exhibits an elbow-like structure around $\underline{\pi_1}$. To estimate $\underline{\pi_1}$, one should select the smallest $\gamma$ for which the distance remains small. This leads to the estimator discussed in \ref{est}. The consistency of the estimator requires a criterion stronger than $||F_n - F ||_{ \infty} \stackrel{a.s.}{\to}  0 $ as $n \to \infty $. Specifically, we need $\sqrt{n} d_n ( F_n , F) $ to be stochastically bounded for the estimator in \ref{est} to be consistent. The following theorem explains the consistency of the estimator in \ref{est}.
\begin{tcolorbox}
\begin{theorem}\label{consistency}
If $ \sqrt{n} d_n ( F_n , F ) $ is stochastically bounded and $ c_n \to \infty $ and $ \lim_{ n \to \infty} \frac{c_n}{ \sqrt{n} } = 0 $ (i.e. $c_n = o ( \sqrt{n} ) $), then $ \hat{\pi}_{1}^{c_n} \stackrel{P}{\to} \text{\underbar{$\pi_1$}} $ as $n \to \infty$.  
\end{theorem} 
\end{tcolorbox}
The proof of theorem \ref{consistency} can be done using the idea similar to \cite{patra2016estimation}. The proof of this theorem is given in appendix. The criterion $\sqrt{n} d_n (F_n , F) $ remains stochastically bounded is stronger than $||F_n - F||_{ \infty} \stackrel{a.s.}{\to} 0 $ as $ n \to \infty  $ which is required by \ref{EMP_CDF}. We highlight a few cases from \cite{billingsley2017probability} where this condition is met.
\subsection{Mixing sequence} 
As per the definitions of \cite{billingsley2013convergence}, consider two measures of dependence between two sigma algebras $\mathscr{F}$ and $\mathscr{G}$. 
\begin{itemize}[label=\ding{70}]
    \item $\alpha_n = \sup \left\{  | P( A \cap B) - P(A) P(B) | \: : \: A \in \mathscr{F}, B \in \mathscr{G}  \right\}$
    \item $\rho_n = \sup \left\{  |E [ \xi \eta ] |  \: : \: E [ \xi ] = E [ \eta ] = 0 , \: || \xi|| \leq 1 , || \eta || \leq 1 , \: \sigma (\xi )  \subseteq \mathscr{F}, \sigma (\eta ) \subseteq \mathscr{G}  \right\}$
\end{itemize} 
It is well known that, $ \alpha_n \leq \rho_n $. 
\begin{itemize}[label=\ding{226}]
    \item A sequence of stationary random variables $X_1 , X_2 , \dots $ is said to be $\alpha$-mixing if for $ \mathscr{F} = \sigma (X_1, \dots , X_k ) $ and $ \mathscr{G} = \sigma (X_{k+n}, X_{ k+n+1} , ... )$, $ \lim_{ n \to \infty} \alpha_n = 0 $. 
    \item A sequence of stationary random variables $X_1 , X_2 , \dots $ is said to be $\rho$-mixing if $ \lim_{ n \to \infty} \rho_n = 0 $. 
\end{itemize} 
The requirement for a mixing sequence is more rigorous than just being stationary and ergodic. For various types of mixing sequences, the criterion $ \sqrt{n} || F_n - F||_{ \infty} = O_P (1) $ is satisfied. We highlight a few examples below. 
\begin{itemize}
    \item For a stationary sequence of random variables $X_1 , X_2 , \dots $, we define a sequence $ \{ \rho_n (x) \}_{ n \geq 1 } $ as follows. 
    $$ \rho_{ | i - j | } (x)  = P( X_i \leq x, X_j \leq x ) - (F(x))^2 $$ \
    If for every $x$, $ \sum\limits_{ n = 1}^{ \infty} \rho_n (x) < \infty $, then we will have $ \sqrt{n} || F_n - F ||_{ \infty} = O_P (1) $ as  a consequence of theorem 19.2 from \cite{billingsley2013convergence}. 
    \item If the stationary sequence of random variables $X_1 , X_2 , ...$ is $ \alpha$-mixing with $ \alpha_n = O( n^{ -5} )$, then $ \sqrt{n} || F_n - F ||_{ \infty} = O_P (1) $. 
\end{itemize} 
Many real-life applications frequently involve examples that meet the assumptions outlined in the previous discussion. For instance, if a stationary sequence is either m-dependent or block dependent, it will satisfy both of the criteria mentioned above. Also, by \cite{billingsley2017probability}, a stationary markov chain with finite state space is $\alpha $-mixing for $ \alpha_n = K \theta^n $ for some $ \theta $ with $|\theta | < 1 $. Thus, it also satisfies the above criterion.
\subsection{Simulation results}
According to Theorem \ref{consistency}, a sequence $\{c_n\}$ is required such that $c_n \to \infty$ as $n \to \infty$ and $\lim_{n \to \infty} \frac{c_n}{\sqrt{n}} = 0$. \cite{patra2016estimation} emphasized the importance of selecting an appropriate $c_n$ for the effective performance of the estimator $\hat{\pi}_{1}^{c_n}$. Their recommendation, based on finite sample performance, suggests $c_n = 0.1 \log (\log n)$ as a useful choice. Additionally, they proposed a cross-validation-based method for choosing the optimal $c_n$. For our simulations, we employed 10-fold cross-validation to determine the optimal $c_n$, and also tested the fixed choice $c_n = 0.1 \log (\log n)$. We have considered some models where $\sqrt{n} \|F_n - F\|_{\infty}$ is stochastically bounded. \\[0.04 in]
Our primary focus is the AR(1) model described in \ref{AR1}. We consider a sample size of \( n = 1000 \) and p-values derived from a mixture of \( N(0,1) \) and \( N(2,1) \) distributions. For the AR(1) parameter \( a \), we examine four values: \( a = 0.1 \), \( 0.2 \), \( 0.5 \), and \( 0.75 \). The performance of the estimator for \( a = 0.1 \) and \( a = 0.2 \) is detailed in Table \ref{AR1_a0102} below. \\[0.04 in]
The estimator \( \hat{\pi}_{1}^{c_n} \) estimates \( \pi_1 = 1 - \pi_0 \). In Table \ref{AR1_a0102}, we present the estimated values of \( \pi_0 \) (i.e., \( 1 - \hat{\pi}_{1}^{c_n} \)), rounded to three decimal places.

\begin{center}
    \begin{table}[ht]
\centering
\caption{AR(1) model with $ a = 0.1$ and $0.2$}
\label{AR1_a0102}
\begin{tabular}{|l|ll|ll|}
\hline
\textbf{} &
  \multicolumn{2}{c|}{\textbf{$ a = 0.1$}} &
  \multicolumn{2}{c|}{\textbf{$ a = 0.2$}} \\ \hline
\textbf{$\pi_0$} &
  \multicolumn{1}{l|}{\textbf{$ \hat{\pi}_0$ (Fixed)}} &
  \textbf{$ \hat{\pi}_0$ (CV)} &
  \multicolumn{1}{l|}{\textbf{$ \hat{\pi}_0$ (Fixed)}} &
  \textbf{$ \hat{\pi}_0$ (CV)} \\ \hline
0.5  & \multicolumn{1}{l|}{0.544} & 0.500 & \multicolumn{1}{l|}{0.547} & 0.512 \\ \hline
0.6  & \multicolumn{1}{l|}{0.636} & 0.596 & \multicolumn{1}{l|}{0.637} & 0.599 \\ \hline
0.65 & \multicolumn{1}{l|}{0.682} & 0.652 & \multicolumn{1}{l|}{0.685} & 0.644 \\ \hline
0.7  & \multicolumn{1}{l|}{0.729} & 0.698 & \multicolumn{1}{l|}{0.724} & 0.700 \\ \hline
0.75 & \multicolumn{1}{l|}{0.770} & 0.729 & \multicolumn{1}{l|}{0.773} & 0.738 \\ \hline
0.8  & \multicolumn{1}{l|}{0.818} & 0.795 & \multicolumn{1}{l|}{0.816} & 0.786 \\ \hline
0.85 & \multicolumn{1}{l|}{0.860} & 0.841 & \multicolumn{1}{l|}{0.859} & 0.835 \\ \hline
0.9  & \multicolumn{1}{l|}{0.899} & 0.886 & \multicolumn{1}{l|}{0.898} & 0.883 \\ \hline
0.95 & \multicolumn{1}{l|}{0.945} & 0.935 & \multicolumn{1}{l|}{0.940} & 0.927 \\ \hline
0.98 & \multicolumn{1}{l|}{0.963} & 0.969 & \multicolumn{1}{l|}{0.966} & 0.953 \\ \hline
\end{tabular}
\end{table}
\end{center} 
For \( a = 0.1 \) and \( a = 0.2 \), there is minimal difference between the estimated values obtained using fixed and cross-validation approaches. The estimators \( \hat{\pi}_0 \) remain close to the true value of \( \pi_0 \). However, as we increase the AR(1) coefficient to \( a = 0.5 \) and \( a = 0.75 \), the gap between the true value of \( \pi_0 \) and \( \hat{\pi}_0 \) widens due to the increased level of dependence. This indicates that the rate of convergence slows as the coefficient \( a \) increases. 
\begin{center}
    \begin{table}[ht]
\centering
\caption{AR(1) model with $ a = 0.5$ and $0.75$}
\label{AR1_a05075}
\begin{tabular}{|l|ll|ll|}
\hline
\multicolumn{1}{|c|}{\textbf{}} &
  \multicolumn{2}{c|}{\textbf{$ a = 0.5$}} &
  \multicolumn{2}{c|}{\textbf{$ a = 0.75$}} \\ \hline
\textbf{$\pi_0$} &
  \multicolumn{1}{l|}{\textbf{$ \hat{\pi}_0$ (Fixed)}} &
  \textbf{$ \hat{\pi}_0$ (CV)} &
  \multicolumn{1}{l|}{\textbf{$ \hat{\pi}_0$ (Fixed)}} &
  \textbf{$ \hat{\pi}_0$ (CV)} \\ \hline
0.5  & \multicolumn{1}{l|}{0.556} & 0.499 & \multicolumn{1}{l|}{0.525} & 0.440 \\ \hline
0.6  & \multicolumn{1}{l|}{0.636} & 0.583 & \multicolumn{1}{l|}{0.564} & 0.474 \\ \hline
0.65 & \multicolumn{1}{l|}{0.668} & 0.629 & \multicolumn{1}{l|}{0.588} & 0.495 \\ \hline
0.7  & \multicolumn{1}{l|}{0.706} & 0.638 & \multicolumn{1}{l|}{0.608} & 0.517 \\ \hline
0.75 & \multicolumn{1}{l|}{0.734} & 0.687 & \multicolumn{1}{l|}{0.627} & 0.548 \\ \hline
0.8  & \multicolumn{1}{l|}{0.772} & 0.706 & \multicolumn{1}{l|}{0.645} & 0.552 \\ \hline
0.85 & \multicolumn{1}{l|}{0.801} & 0.756 & \multicolumn{1}{l|}{0.664} & 0.558 \\ \hline
0.9  & \multicolumn{1}{l|}{0.837} & 0.782 & \multicolumn{1}{l|}{0.684} & 0.593 \\ \hline
0.95 & \multicolumn{1}{l|}{0.866} & 0.806 & \multicolumn{1}{l|}{0.701} & 0.629 \\ \hline
0.98 & \multicolumn{1}{l|}{0.885} & 0.839 & \multicolumn{1}{l|}{0.714} & 0.630 \\ \hline
\end{tabular}
\end{table}
\end{center}
It is interesting to observe that the discrepancy between \( \hat{\pi}_0 \) and \( \pi_0 \) becomes more pronounced as \( \pi_0 \) increases. The estimator's performance is satisfactory for smaller values of \( \pi_0 \) (up to \( \pi_0 = 0.6 \)). Moreover, using a cross-validation-based estimator provides little or no improvement in this range. In these cases, the fixed choice of \( c_n = 0.1 \log (\log n) \) also yields good results. \\[0.04 in] 
We now consider the m-dependent scenario. As before, given a sample \( \{ \xi_i \}_{i \geq 1} \), observations are generated according to \( X_i \propto \sum\limits_{j = i}^{i+m} \left( (m+1) - (j-i) \right) \xi_j \) and are standardized to have unit variance. A mixture of \( N(0,1) \) and \( N(2,1) \) variables is then formed from this sequence \( \{ X_i \}_{i \geq 1} \), and the p-values are computed for a one-sided test of the mean. We conduct simulations with \( n = 1000 \) and consider \( m = 1, 2, \) and \( 5 \). The estimated values of \( \pi_0 \) are shown in the table \ref{m-dep}. 
\begin{center}
    \begin{table}[ht]
\centering
\caption{Estimator of null proportion for m-dependent model with m = 1,2,5}
\label{m-dep}
\begin{tabular}{|l|ll|ll|ll|}
\hline
\textbf{} & \multicolumn{2}{c|}{\textbf{m=1}}  & \multicolumn{2}{c|}{\textbf{m=2}}  & \multicolumn{2}{c|}{\textbf{m=5}}  \\ \hline
\textbf{$\pi_0$} &
  \multicolumn{1}{l|}{\textbf{$\hat{\pi}_0$ ( Fixed)}} &
  \textbf{$\hat{\pi}_0$ ( CV)} &
  \multicolumn{1}{l|}{\textbf{$\hat{\pi}_0$ ( Fixed)}} &
  \textbf{$\hat{\pi}_0$ ( CV)} &
  \multicolumn{1}{l|}{\textbf{$\hat{\pi}_0$ ( Fixed)}} &
  \textbf{$\hat{\pi}_0$ ( CV)} \\ \hline
0.5       & \multicolumn{1}{l|}{0.540} & 0.460 & \multicolumn{1}{l|}{0.512} & 0.425 & \multicolumn{1}{l|}{0.449} & 0.345 \\ \hline
0.6       & \multicolumn{1}{l|}{0.599} & 0.508 & \multicolumn{1}{l|}{0.550} & 0.459 & \multicolumn{1}{l|}{0.466} & 0.368 \\ \hline
0.65      & \multicolumn{1}{l|}{0.622} & 0.549 & \multicolumn{1}{l|}{0.568} & 0.470 & \multicolumn{1}{l|}{0.474} & 0.388 \\ \hline
0.7       & \multicolumn{1}{l|}{0.649} & 0.571 & \multicolumn{1}{l|}{0.588} & 0.503 & \multicolumn{1}{l|}{0.482} & 0.408 \\ \hline
0.75      & \multicolumn{1}{l|}{0.674} & 0.602 & \multicolumn{1}{l|}{0.605} & 0.518 & \multicolumn{1}{l|}{0.491} & 0.378 \\ \hline
0.8       & \multicolumn{1}{l|}{0.701} & 0.627 & \multicolumn{1}{l|}{0.622} & 0.543 & \multicolumn{1}{l|}{0.497} & 0.418 \\ \hline
0.85      & \multicolumn{1}{l|}{0.728} & 0.679 & \multicolumn{1}{l|}{0.638} & 0.543 & \multicolumn{1}{l|}{0.505} & 0.410 \\ \hline
0.9       & \multicolumn{1}{l|}{0.750} & 0.688 & \multicolumn{1}{l|}{0.656} & 0.549 & \multicolumn{1}{l|}{0.513} & 0.423 \\ \hline
0.95      & \multicolumn{1}{l|}{0.773} & 0.703 & \multicolumn{1}{l|}{0.672} & 0.580 & \multicolumn{1}{l|}{0.522} & 0.430 \\ \hline
0.98      & \multicolumn{1}{l|}{0.786} & 0.718 & \multicolumn{1}{l|}{0.682} & 0.594 & \multicolumn{1}{l|}{0.525} & 0.437 \\ \hline
\end{tabular}
\end{table}
\end{center}
Due to the very slow rate of convergence, the deviation between the true and estimated values of \( \pi_0 \) is evident in this setup. This discrepancy becomes even more pronounced as the degree of dependence increases, which occurs when the value of \( m \) is raised. Additionally, in the m-dependent case where equal weights are assigned to neighbors (i.e., \( X_i = \frac{1}{\sqrt{m+1}} \sum\limits_{j=i}^{i+m} \xi_j \)), the deviation between the true and estimated values is even more substantial. This is clearly illustrated in the table \ref{m-dep_equal}. It is noteworthy that, due to the increased dependence among observations, the performance of the cross-validation-based choice of $c_n$ is inferior to the fixed choice regarding the deviation between the true and estimated values of $\pi_0$. The assumption of independence is crucial for the validity of cross-validation-based methods, which is why the fixed choice of $c_n$ performs better in scenarios with a higher degree of dependence. Consequently, the extensive time spent on cross-validation appears to be unproductive in this context.
\begin{center}

\begin{table}[ht]
\centering
\caption{Equally weighted m-dependent model with m = 1,2,5}
\label{m-dep_equal}
\begin{tabular}{|l|ll|ll|ll|}
\hline
     & \multicolumn{2}{c|}{\textbf{m=1}}  & \multicolumn{2}{c|}{\textbf{m=2}}  & \multicolumn{2}{c|}{\textbf{m=5}}  \\ \hline
\textbf{$\pi_0$} &
  \multicolumn{1}{l|}{\textbf{$\hat{\pi}_0$ (Fixed)}} &
  \textbf{$\hat{\pi}_0$ (CV)} &
  \multicolumn{1}{l|}{\textbf{$\hat{\pi}_0$ (Fixed)}} &
  \textbf{$\hat{\pi}_0$ (CV)} &
  \multicolumn{1}{l|}{\textbf{$\hat{\pi}_0$ (Fixed)}} &
  \textbf{$\hat{\pi}_0$ (CV)} \\ \hline
0.5  & \multicolumn{1}{l|}{0.443} & 0.381 & \multicolumn{1}{l|}{0.449} & 0.398 & \multicolumn{1}{l|}{0.444} & 0.386 \\ \hline
0.6  & \multicolumn{1}{l|}{0.501} & 0.448 & \multicolumn{1}{l|}{0.498} & 0.436 & \multicolumn{1}{l|}{0.496} & 0.429 \\ \hline
0.65 & \multicolumn{1}{l|}{0.516} & 0.482 & \multicolumn{1}{l|}{0.522} & 0.472 & \multicolumn{1}{l|}{0.515} & 0.463 \\ \hline
0.7  & \multicolumn{1}{l|}{0.550} & 0.495 & \multicolumn{1}{l|}{0.546} & 0.499 & \multicolumn{1}{l|}{0.537} & 0.467 \\ \hline
0.75 & \multicolumn{1}{l|}{0.564} & 0.528 & \multicolumn{1}{l|}{0.570} & 0.523 & \multicolumn{1}{l|}{0.561} & 0.513 \\ \hline
0.8  & \multicolumn{1}{l|}{0.588} & 0.546 & \multicolumn{1}{l|}{0.585} & 0.553 & \multicolumn{1}{l|}{0.589} & 0.525 \\ \hline
0.85 & \multicolumn{1}{l|}{0.608} & 0.572 & \multicolumn{1}{l|}{0.611} & 0.563 & \multicolumn{1}{l|}{0.608} & 0.538 \\ \hline
0.9  & \multicolumn{1}{l|}{0.632} & 0.591 & \multicolumn{1}{l|}{0.627} & 0.576 & \multicolumn{1}{l|}{0.608} & 0.581 \\ \hline
0.95 & \multicolumn{1}{l|}{0.648} & 0.614 & \multicolumn{1}{l|}{0.640} & 0.595 & \multicolumn{1}{l|}{0.636} & 0.580 \\ \hline
0.98 & \multicolumn{1}{l|}{0.658} & 0.621 & \multicolumn{1}{l|}{0.652} & 0.615 & \multicolumn{1}{l|}{0.649} & 0.597 \\ \hline
\end{tabular}
\end{table}
\end{center} 
We now shift our focus to the block-dependent setup where the requirement of stationarity is not satisfied. In this context, the estimator $\hat{\pi}_{1}^{c_n}$ appears to perform better compared to the m-dependent case. We fixed $n=1000$ and examined two scenarios: one with 10 blocks of size 100 and another with 200 blocks of size 5. Clearly, the latter scenario is closer to independence than the former. The estimators are presented in Table \ref{Block_dep_PS}. As before, we considered 10 different correlation values for the block-dependent case, namely $0.1, \: 0.2, \: \dots, \: 0.9, \: 0.95$.
\begin{center}
    \begin{table}[ht]
\centering
\caption{Block dependent case with $n=1000$ }
\label{Block_dep_PS}
\begin{tabular}{|l|ll|ll|}
\hline
\textbf{} &
  \multicolumn{2}{c|}{\textbf{Block of size 5}} &
  \multicolumn{2}{c|}{\textbf{Block of size 100}} \\ \hline
\textbf{$\pi_0$} &
  \multicolumn{1}{l|}{\textbf{$\hat{\pi}_0$ (Fixed)}} &
  \textbf{$\hat{\pi}_0$ (CV)} &
  \multicolumn{1}{l|}{\textbf{$\hat{\pi}_0$ (Fixed)}} &
  \textbf{$\hat{\pi}_0$ (CV)} \\ \hline
0.5  & \multicolumn{1}{l|}{0.531} & 0.481 & \multicolumn{1}{l|}{0.540} & 0.488 \\ \hline
0.6  & \multicolumn{1}{l|}{0.616} & 0.568 & \multicolumn{1}{l|}{0.628} & 0.578 \\ \hline
0.65 & \multicolumn{1}{l|}{0.657} & 0.621 & \multicolumn{1}{l|}{0.671} & 0.635 \\ \hline
0.7  & \multicolumn{1}{l|}{0.696} & 0.671 & \multicolumn{1}{l|}{0.716} & 0.670 \\ \hline
0.75 & \multicolumn{1}{l|}{0.737} & 0.697 & \multicolumn{1}{l|}{0.755} & 0.690 \\ \hline
0.8  & \multicolumn{1}{l|}{0.773} & 0.753 & \multicolumn{1}{l|}{0.799} & 0.749 \\ \hline
0.85 & \multicolumn{1}{l|}{0.815} & 0.774 & \multicolumn{1}{l|}{0.836} & 0.807 \\ \hline
0.9  & \multicolumn{1}{l|}{0.851} & 0.817 & \multicolumn{1}{l|}{0.876} & 0.849 \\ \hline
0.95 & \multicolumn{1}{l|}{0.881} & 0.877 & \multicolumn{1}{l|}{0.911} & 0.890 \\ \hline
0.98 & \multicolumn{1}{l|}{0.896} & 0.900 & \multicolumn{1}{l|}{0.920} & 0.892 \\ \hline
\end{tabular}
\end{table}
\end{center}
Interestingly, the performance of the cross-validation-based estimator deteriorates when moving from $200$ blocks to $10$ blocks, due to the increased degree of dependence. In contrast, the fixed choice of $c_n$ provides a relatively satisfactory estimate of $\pi_0$. However, using these estimates may not ensure conservative control over the false discovery rate (FDR), as they tend to underestimate $\pi_0$ at higher values.
\section{Remarks}
The simulation results indicate that both Storey's estimator and the estimator proposed by \cite{patra2016estimation} are sensitive to the dependence structure in the data. Storey's estimator is known to be a biased estimator of \( \pi_0 \), whereas the estimator \( \hat{\pi}_{1}^{c_n} \) is consistent under conditions of weak dependence. However, its rate of convergence is notably slow, and this rate decreases even further in the presence of dependence. Even relatively simple dependence structures, such as m-dependence, can significantly distort the performance of these estimators. Therefore, further research is needed to develop a consistent estimator of the proportion of non-null hypotheses that remains reliable under weak dependence.

\section{Appendix}
\subsection{Proof of theorem \ref{th_cons}}
\textbf{Case 1}: \{$\gamma \geq\text{  \underbar{$\pi_1$}} $\}  \\[0.05 in] 
Let, \begin{equation}\label{F1gamma}F_{1,\gamma} = \frac{ F - (1- \gamma) F_0 }{ \gamma} \end{equation}
Then, $\gamma d_n ( \hat{F}_{1,n}, \hat{F}_{1,\gamma} ) = d_n ( F_n , F) $. \\[0.08 in]  
Consider the set $ \Gamma = \Big\{ \gamma \in [0,1] \: | \: \{ \frac{ F - ( 1- \gamma) F_0 }{ \gamma } \Big\}\text{ is a valid cdf } \}$. We'll first show that, $ \Gamma $ is a convex set. \\[0.05 in]
Suppose $ \gamma_1 , \gamma_2 \in \Gamma $ so that, $ \frac{ F - (1 - \gamma_i) F_0}{\gamma_i} $ is a valid cdf for $i=1,2$. \\[0.05 in] 
For any $ 0 \leq \alpha \leq 1 $, \\[0.08 in]
$ \frac{ F - ( 1 - \alpha \gamma_1 - (1- \alpha) \gamma_2)  F_0}{ \alpha \gamma_1 + (1- \alpha) \gamma_2 } $ \\[0.05 in]
=$\Big( \frac{ \alpha \gamma_1 }{\alpha \gamma_1 + (1- \alpha) \gamma_2} \Big) \left(  \frac{ F - (1- \gamma_1) F_0 }{ \gamma_1 } \right) + 
\Big( \frac{ (1-\alpha) \gamma_2 }{\alpha \gamma_1 + (1- \alpha) \gamma_2} \Big) \left(  \frac{ F - (1- \gamma_2) F_0 }{ \gamma_1 } \right)$ \\[0.08 in]
Since convex combination of cdfs should be cdf, we can now say that, $\frac{ F - ( 1 - \alpha \gamma_1 - (1- \alpha) \gamma_2)  F_0}{ \alpha \gamma_1 + (1- \alpha) \gamma_2 } $ is a valid cdf and hence $ \alpha \gamma_1 + (1- \alpha) \gamma_2 \in \Gamma $. \\[0.05 in]
This implies, for any $ \text{  \underbar{$\pi_1$}} \leq \gamma \leq  1 $, $ \frac{ F - ( 1- \gamma) F_0}{ \gamma} $ is a valid cdf. Thus, $ F_{ 1 , \gamma} $ defined in \ref{F1gamma} is a valid cdf. \\[0.05 in] 
So, by the definition of $\hat{F}_{2,n} $, it can be said that, 
$$
\gamma d_n ( \hat{F}_{1,n}, \hat{F}_{2,n} ) \leq \gamma d_n ( \hat{F}_{1,n}, \hat{F}_{1,\gamma} ) =d_n ( F_n , F) \quad \forall \: \gamma \geq \text{  \underbar{$\pi_1$}}
$$
Since, $ d_n (F_n , F) \leq ||F_n  - F||_{ \infty} $, $ ||F_n - F||_{ \infty} \stackrel{a.s.}{\to}  0 $ as $ n \to \infty $ implies, $ d_n (F_n , F) \stackrel{a.s.}{\to}  0 $ as $ n \to \infty $. Hence, 
$$ \gamma d_n ( \hat{F}_{1,n}, \hat{F}_{2,n} ) \stackrel{a.s.}{\to}  0 \text{ as } n \to \infty $$
\textbf{Case 2}: \{$\gamma <\text{  \underbar{$\pi_1$}} $\}  \\[0.05 in] 
$\hat{F_{1,n}} $ is not a valid cdf. by the definition of $\text{  \underbar{$\pi_1$}}$. Also, if $ || F_n - F||_{ \infty} \stackrel{a.s.}{\to}  0 $, then $ \hat{ F_{1,n}} \stackrel{a.s.}{\to}  F_{1, \gamma} $ as $ n \to \infty$. \\[0.05 in] 
So, for sufficientlly large $n$, the difference between $ \hat{F}_{1,n} $ and $\hat{F}_{2,n}$ is substantially high as $ \hat{F}_{1,n} $ converges to $ F_{1 , \gamma} $ which is not a valid cdf and by definition, $\hat{F}_{2,n} $ is a cdf. \\[0.05 in]
Combining the two observations in case 1 and case 2, we can finally conclude that, 
   \begin{equation*}
\gamma d_n ( \hat{F}_{1,n}, \hat{F}_{2,n} ) \stackrel{a.s.}{\to}  
\begin{cases} 
0 \qquad \text{  if } \gamma \geq\text{  \underbar{$\pi_1$}} \\
> 0, \: \:  \text{  if } \gamma <\text{  \underbar{$\pi_1$}}
\end{cases} 
\end{equation*}

\subsection{Proof of theorem \ref{consistency}} 
Fix any $ \epsilon > 0 $. We need to show that, $ \lim_{ n \to \infty} P( | \hat{\pi}_{1}^{c_n} -\text{  \underbar{$\pi_1$}} | \geq \epsilon ) = 0 $ \\[0.05 in]
We will first prove that, $ \lim_{ n \to \infty} P(  \hat{\pi}_{1}^{c_n} -\text{  \underbar{$\pi_1$}} \leq - \epsilon ) = 0 $. \\[0.05 in] 
If $\text{  \underbar{$\pi_1$}} \leq \epsilon$, then $ \text{  \underbar{$\pi_1$}} - \epsilon \leq 0 $ and hence it is obvious. So, assume that $\text{  \underbar{$\pi_1$}}> \epsilon$. \\[0.05 in]
Recall that, $$\text{  \underbar{$\pi_1$}} = \inf \: \{ \: \gamma \in (0, 1 ] \: : \: \Big\{ \frac{ F - (1- \gamma) F_0 }{ \gamma } \Big\} \text{ is a CDF} \}$$ 
So, for $ \gamma <\text{  \underbar{$\pi_1$}}, \: \: \frac{ F - (1- \gamma) F_0 }{ \gamma }$ is not a cdf and for any such $\gamma$, $\hat{F}_{1,n}, \hat{F}_{2,n}$ are quite different.\\[0.05 in]
Let, $ \hat{F}_{1,n} ( \gamma), \hat{F}_{2,n} ( \gamma) $ be the value of the estimators for a fixed $\gamma$.\\[0.05 in]
Consider $\gamma = (\text{  \underbar{$\pi_1$}} - \epsilon) $, and we can conclude that, $ (\text{  \underbar{$\pi_1$}} - \epsilon) d_n \left( \hat{F}_{1,n} (\text{  \underbar{$\pi_1$}} - \epsilon) , \hat{F}_{2,n} ( \text{  \underbar{$\pi_1$}} - \epsilon) \right)$ is bounded away from zero. \\[0.05 in] 
Recall that, 
$$ \hat{\pi}_{1}^{c_n} = \inf \left\{  \gamma \in (0,1] \: : \: \gamma d_n ( \hat{F}_{1,n}, \hat{F}_{2,n} ) \leq \frac{c_n}{\sqrt{n}} \right\}$$ 
Consider the set $ A_n = \left\{  \gamma \in (0,1] \: : \: \gamma d_n ( \hat{F}_{1,n}, \hat{F}_{2,n} ) \leq \frac{c_n}{\sqrt{n}} \right\}$. Clearly $A_n $ is a convex set. \\[0.05 in] 
So, $ P( \hat{\pi}_{1}^{c_n}  \leq\text{  \underbar{$\pi_1$}} - \epsilon ) 
\leq P \Big( (\text{  \underbar{$\pi_1$}}- \epsilon)d_n \left( \hat{F}_{1,n} (\text{  \underbar{$\pi_1$}}- \epsilon), \hat{F}_{2,n} (\text{  \underbar{$\pi_1$}}- \epsilon \right) \leq \frac{ c_n}{ \sqrt{n}} \Big) $\\[0.05 in] 
Since, $(\text{  \underbar{$\pi_1$}} - \epsilon) d_n \left( \hat{F}_{1,n} (\text{  \underbar{$\pi_1$}} - \epsilon) , \hat{F}_{2,n} (\text{  \underbar{$\pi_1$}} - \epsilon) \right) $ is bounded away from zero and $\lim_{ n \to \infty} \frac{c_n}{ \sqrt{n}} = 0 $, we can say that, 
$$ \lim_{ n \to \infty } P \Big( (\text{  \underbar{$\pi_1$}}- \epsilon)d_n \left( \hat{F}_{1,n} (\text{  \underbar{$\pi_1$}}- \epsilon), \hat{F}_{2,n} (\text{  \underbar{$\pi_1$}}- \epsilon) \right) \leq \frac{ c_n}{ \sqrt{n}} \Big)= 0 $$ 
Thus, $ \lim_{ n \to \infty} P( \hat{\pi}_{1}^{c_n}  \leq\text{  \underbar{$\pi_1$}} - \epsilon ) = 0  $ \\[0.05 in] 
On the other hand, $ \{ \hat{\pi}_{1}^{c_n}  -\text{  \underbar{$\pi_1$}} > \epsilon )\}$ implies, $ \Big\{ (\text{  \underbar{$\pi_1$}}+ \epsilon )d_n \left(\hat{F}_{1,n} ((\text{  \underbar{$\pi_1$}}+ \epsilon)), \hat{F}_{2,n} ((\text{  \underbar{$\pi_1$}}+ \epsilon) \right) > \frac{c_n}{ \sqrt{n}} \Big\}$. So, \\[0.05 in ] 
$ P( \hat{\pi}_{1}^{c_n}  -\text{  \underbar{$\pi_1$}} > \epsilon ) \leq 
P \Big(  \sqrt{n} d_n \left(\hat{F}_{1,n} ((\text{  \underbar{$\pi_1$}}+ \epsilon)), \hat{F}_{2,n} ((\text{  \underbar{$\pi_1$}}+ \epsilon) \right) > \frac{c_n}{(\text{  \underbar{$\pi_1$}}+ \epsilon )} \Big) $  \\[0.08 in] 
From the proof of theorem \ref{th_cons}, for any $ \gamma \in (\text{  \underbar{$\pi_1$}},1], \: \: \gamma d_n ( \hat{F}_{1,n}, \hat{F}_{2,n} ) \leq d_n ( F_n , F) $. \\[0.08 in] 
So, $ 
 P( \hat{\pi}_{1}^{c_n}  -\text{  \underbar{$\pi_1$}} > \epsilon ) \leq P \Big( \sqrt{n} d_n  ( F_n , F) > \frac{c_n}{(\text{  \underbar{$\pi_1$}}+ \epsilon )} \Big)   $  \\[0.07 in]  
Since, $ \lim_{ n \to \infty } c_n = \infty $ and $ \sqrt{n} d_n ( F_n , F) $ is stochastically bounded, we can conclude that, 
$$ \lim_{ n \to \infty } P \left( \sqrt{n} d_n  ( F_n , F) > \frac{c_n}{(\text{  \underbar{$\pi_1$}}+ \epsilon )} \right) = 0 $$ 
Thus, $ \lim_{ n \to \infty } P( \hat{\pi}_{1}^{c_n}  -\text{  \underbar{$\pi_1$}} > \epsilon ) = 0 $ and this completes the proof. 

\bibstyle{plain}
\bibliography{references.bib}

\end{document}